\documentclass[11pt,article]{article}

\usepackage{amsmath}
\usepackage{amsthm}
  \usepackage{paralist}
  \usepackage{graphics} 
  \usepackage{epsfig} 
\usepackage{graphicx}
\usepackage{caption}
\usepackage{subcaption}
\usepackage{epstopdf}
 \usepackage[colorlinks=true]{hyperref}
 \usepackage{multirow}
 \usepackage[font=scriptsize,labelfont=sf]{caption}
\usepackage[margin=2cm]{caption}
\input{amssym.tex}

\newtheorem{theorem}{Theorem}[section]
\newtheorem{corollary}{Corollary}
\newtheorem{lemma}[theorem]{Lemma}
\newtheorem{proposition}{Proposition}[section]

 \numberwithin{equation}{section}
\newtheorem{remark}{Remark}[section]

\newcommand{\keywords}

\def\bc{\begin{center}}       \def\ec{\end{center}}
\def\ba{\begin{array}}        \def\ea{\end{array}}
\def\be{\begin{equation}}     \def\ee{\end{equation}}
\def\bea{\begin{eqnarray}}    \def\eea{\end{eqnarray}}
\def\beaa{\begin{eqnarray*}}  \def\eeaa{\end{eqnarray*}}

\def\mathbb{\Bbb}

\begin{document}

\title{\bf Time--periodic and stable patterns of a two--competing--species Keller--Segel chemotaxis model: effect of cellular growth\thanks{\emph{Discrete Contin. Dyn. Syst. Ser. B}, 22 (2017), no. 9, 3547-3574. }}
\author{Qi Wang \thanks{Department of Mathematics, Southwestern University of Finance and Economics, Chengdu, Sichuan 611130, China.  Email: {\tt qwang@swufe.edu.cn}.  Corresponding author.  QW is partially supported by NSF-China (Grant 11501460)}, Jingyue Yang \thanks{School of Finance, Southwestern University of Finance and Economics, Chengdu, Sichuan 611130, China.  Email: {\tt yjy@2011.swufe.edu.cn}.}, Lu Zhang \thanks{Department of Mathematics, Southern Methodist University, Dallas, Texas 75275, USA.  Email: {\tt luzhang@smu.edu.cn}.}\\
}
\date{}
\maketitle


\abstract
This paper investigates the formation of time--periodic and stable patterns of a two--competing--species Keller--Segel chemotaxis model with a focus on the effect of cellular growth.  We carry out rigorous Hopf bifurcation analysis to obtain the bifurcation values, spatial profiles and time period associated with these oscillating patterns.  Moreover, the stability of the periodic solutions is investigated and it provides a selection mechanism of stable time--periodic mode which suggests that only large domains support the formation of these periodic patterns.  Another main result of this paper reveals that cellular growth is responsible for the emergence and stabilization of the oscillating patterns observed in the $3\times3$ system, while the system admits a Lyapunov functional in the absence of cellular growth.  Global existence and boundedness of the system in 2D are proved thanks to this Lyapunov functional.  Finally, we provide some numerical simulations to illustrate and support our theoretical findings.

\textbf{Keywords: Two species chemotaxis model, time--periodic pattern, stable pattern, Hopf bifurcation, stability analysis, Lyapunov functional.}

\section{Introduction and preliminary results}
In this paper, we continue our study in \cite{WZYH} of the following system of $(u,v,w)=(u(x,t),v(x,t),w(x,t))$
\begin{equation}\label{11}
\left\{
\begin{array}{ll}
u_t=\nabla \cdot (d_1 \nabla u-\chi u \nabla w)+\mu_1(1-u-a_1v)u,&x \in \Omega,t>0, \\
v_t=\nabla \cdot(d_2\nabla v-\xi v \nabla w)+\mu_2(1-a_2u-v)v,&x \in \Omega,t>0, \\
\tau w_t=\Delta w-\lambda w +u+v,&x \in \Omega,t>0,\\
\frac{\partial u}{\partial \textbf{n}}=\frac{\partial v}{\partial \textbf{n}}=\frac{\partial w}{\partial \textbf{n}}=0,&x\in\partial \Omega,t>0,\\
u(x,0)=u_0(x),v(x,0)=v_0(x),w(x,0)=w_0(x), &x\in \Omega,
\end{array}
\right.
\end{equation}
where $d_i$, $\mu_i$, $a_i$, $i=1,2$, and $\lambda$ are positive constants and $\tau$ is a nonnegative constant; $\Omega$ is a bounded domain in $\mathbb{R}^N$, $N\geq1$ with smooth boundary $\partial \Omega$.  (\ref{11}) is a Keller--Segel type model of chemotaxis, the oriented movement of cellular organisms towards the high concentration region of a chemical released by the cells, and it was proposed by Tello and Winkler in \cite{TW} to study the population dynamics of two competitive biological species attracted by the same nutrition subject to Lotka--Volterra dynamics.  Here, $u(x,t)$ and $v(x,t)$ represent population densities of the two competing species at space--time location $(x,t)\in \Omega\times \mathbb{R}^+$, $w(x,t)$ concentration of the attracting chemical.  It is assumed that both species direct their movement chemotactically along the gradient of chemical concentration over the habitat, hence both $\chi$ and $\xi$ are selected to be positive constants.  Biologically, $\chi$ and $\xi$ measure the strength of chemical attraction to species $u$ and $v$ respectively.  Kinetics of the species are assumed to be of the classical Lotka--Volterra type, hence $a_i$, $i=1,2$, interpret the levels of inter--specific competition and $\mu_i$, $i=1,2$, measure the intrinsic cellular growth.  The chemical is produced by both species at the same rate with no saturation and is consumed by certain enzyme at the rate of $\lambda$ meanwhile.

One of the most interesting phenomena in chemotaxis is the cellular aggregation during which initially homogeneously distributed cells can aggregate and develop into a fruiting body over time.  For example, during the first phase of its developmental cycle, \emph{Dictyostelium discoideum} exists as a single amoeboid cell, then it differentiates into multiple cells which eventually aggregate into a multicellular organism after the growth phase.  In \emph{Dictyostelium} chemotaxis, it was discovered that the aggregating cells of \emph{D. discoideum} are attracted by a chemical called \emph{cyclic AMP} (cAMP), which is synthesized and released by the cells periodically.  See the review papers \cite{Ge,HD}.  It is of great interest to both biologists and mathematicians to understand the initiation and formation of the self--organized oscillating patterns.  Mathematical modeling of chemotaxis dates back to the pioneering works of Patlak \cite{Pt} and Keller--Segel \cite{KS,KS1,KS2}, where a group of parabolic reaction--diffusion systems have been proposed to describe the spatial--temporal behaviors of cellular distribution and chemical concentration.  Diffusion interprets the random movements of the cell and chemical, advection the cellular chemotactic movement, and kinetics the cellular birth--death and chemical degradation--creation.  Mathematical modeling and analysis of chemotaxis have developed substantially since the appearance of works of Keller and Segel.  See the survey papers \cite{HP,Ho,Ho2,Ho3,Wz} and the references cited therein.  We would like to point out that most papers in the literature focus on the studies of chemotaxis model with single bacteria and one chemical stimulus.

Let us now present a brief review of the studies of (\ref{11}) and propose our motivation for the current work.  If $0\leq a_1<1,0\leq a_2<1$, system (\ref{11}) has a unique positive constant steady state
\begin{equation}\label{12}
(\bar u,\bar v,\bar w)=\Big( \frac{1-a_1}{1-a_1a_2}, \frac{1-a_2}{1-a_1a_2},\frac{2-a_1-a_2}{\lambda (1-a_1a_2)} \Big).
\end{equation}
Tello and Winkler \cite{TW} showed that when $\tau=0$, $d_1=d_2=1$ and
\begin{equation}\label{13}
2(\chi+\xi)+a_1\mu_2<\mu_1,~2(\chi+\xi)+a_2\mu_1<\mu_2,
\end{equation}
the constant equilibrium $(\bar u,\bar v,\bar w)$ given by (\ref{12}) is a global attractor of the parabolic--parabolic--elliptic system of (\ref{11}) for any positive initial data $(u_0,v_0,w_0)\in C^0(\bar \Omega)\times C^0(\bar \Omega)\times W^{1,p}(\bar\Omega), p>N$.  Recently, Stinner \emph{et al}. \cite{STW} studied the competitive exclusion of (\ref{11}) with $\tau=0$, under more complicated smallness assumptions on the chemotaxis rates.  In particular, (\ref{11}) has no nonconstant positive steady state when the cellular chemotaxis sensitivity is moderate.

To see that chemotaxis is responsible for the formation of nontrivial patterns in (\ref{11}), we consider the global dynamics of (\ref{11}) with $\chi=\xi=0$, i.e., the following system
\begin{equation}\label{14}
\left\{
\begin{array}{ll}
u_t= d_1 \Delta u +\mu_1(1-u-a_1v)u,&x \in \Omega,t>0, \\
v_t= d_2 \Delta v+\mu_2(1-a_2u-v)v,&x \in \Omega,t>0, \\
\tau w_t=\Delta w-\lambda w +u+v,&x \in \Omega,t>0,\\
\frac{\partial u}{\partial \textbf{n}}=\frac{\partial v}{\partial \textbf{n}}=\frac{\partial w}{\partial \textbf{n}}=0,&x\in\partial \Omega,t>0,\\
u(x,0)=u_0(x),v(x,0)=v_0(x),w(x,0)=w_0(x), &x\in \Omega.
\end{array}
\right.
\end{equation}
It is well known that the first two equations in (\ref{14}) correspond to the weak competition case of classical Lotka--Volterra system for which $(\bar u,\bar v)$ is globally asymptotically stable, then by applying the comparison principle to the $w$--equation we can easily show that $(\bar u,\bar v, \bar w)$ is the global attractor thanks to the classical results in \cite{KiW,LN,WGY}.  Therefore we have the following result.
\begin{proposition}\label{proposition11}
Let $\Omega\subset \mathbb R^N$, $N\geq1$ be a bounded domain with smooth boundary.  Assume that all parameters in (\ref{14}) are positive.  Then for any $(u_0,v_0,w_0)\in C^0(\Omega)\times C^0(\Omega)\times W^{1,p}(\Omega)$, $p>N$, (\ref{14}) has a unique positive classical solution $(u(x,t),v(x,t),w(x,t))$ which is uniformly bounded in time.  Moreover, $(\bar u,\bar v,\bar w)$ is a global attractor of (\ref{14}) such that $\lim_{t\rightarrow \infty} \Vert u(\cdot,t)-\bar u\Vert_{L^\infty}+\Vert v(\cdot,t)-\bar v\Vert_{L^\infty}+\Vert w(\cdot,t)-\bar w\Vert_{L^\infty}=0$.
\end{proposition}
Proposition \ref{proposition11} complements the results in \cite{TW} which requires $a_2<\frac{\mu_2}{\mu_1}<\frac{1}{a_1}$ when $\chi=\xi=0$.  It also suggests that (\ref{14}) does not admit Turing's instability or the so--called diffusion driven instability in the sense that the uniform equilibrium, stable for the ODE system, becomes unstable as a solution to the full reaction--diffusion system (\ref{11}), therefore chemotaxis is responsible for the existence of nonconstant positive solutions for (\ref{11}).  We want to mention that positive steady states of chemotaxis systems with concentrating property are usually adopted to model the cellular aggregation phenomenon.

From the viewpoint of linearized stability analysis, chemotaxis destabilizes spatially homogeneous solutions for reaction--diffusion systems, while diffusion stabilizes homogeneous solutions (unless Turing's instability occurs), therefore one can expect the emergence of spatially inhomogeneous when chemotaxis rate is large.  In \cite{WZYH}, the authors obtained the existence and stability of nontrivial positive steady states to (\ref{11}) over $(0,L)$ through rigorous bifurcation analysis; a selection mechanism of stable wave mode has been proposed to predict the spatial profile of the stable stationary patterns.  Numerical simulations in \cite{WZYH} verify the theoretical findings there and suggest that system (\ref{11}) over $(0,L)$ also admits time--periodic spatial patterns for properly chosen parameters.  It is surmised in \cite{WZYH} that stable oscillating patterns emerge since $(\bar u,\bar v,\bar w)$ loses its stability through Hopf bifurcation.  One of the goals of our current work is to rigorously investigate the formation of stable time--periodic patterns.  In particular, we show that chemotaxis and cellular kinetics are responsible for the formation of temporal oscillating patterns.  We want to mention that the phenomenon of time--periodic oscillations is important not only for reaction--diffusion models in biological and ecological systems (\cite{Dancer} e.g), but almost all other dynamical systems of scientific disciplines such as fluid mechanics \cite{Iu,Iu1,Iu2, Sat1}, lasers \cite{HL} etc.  For example, for the effect of chemotaxis on bacterial strategies, see the discussions in \cite{STW} and the references cited therein.

The rest of this paper is organized as follows.  In Section \ref{section2}, we analyze the linearized stability of $(\bar u,\bar v,\bar w)$ in terms of chemotaxis rate $\chi$.  It is shown that this homogeneous solution loses its stability as $\chi$ surpasses a threshold value $\chi_0$, the minimum of bifurcation values $\chi_k^S$ and $\chi_k^H$ over $\mathbb N^+$, where $\chi_k^S$ and $\chi_k^H$ are chosen such that the stability matrix  (\ref{23}) of $(\bar u,\bar v,\bar w)$ has zero or purely imaginary eigenvalues, respectively.  Section \ref{section3} is devoted to the rigorous Hopf bifurcation analysis of (\ref{11}) over $(0,L)$, for which time--periodic spatial patterns are established.  Our existence results employ Hopf bifurcation theorem for parabolic systems from \cite{A,CR4} etc.  We also investigate the stability of these periodic solutions and establish a selection mechanism for stable oscillating patterns in terms of system parameters.  In Section \ref{section4}, we study the effect of cellular growth on the spatial--temporal dynamics of (\ref{11}).  In particular, we find that when $\mu_1=\mu_2=0$, (\ref{11}) does not admit time--periodic patterns.  Our argument is based on the construction of a time--monotone Lyapunov functional to (\ref{11}).  Moreover, global existence of this problem over 2D is also established provided that the initial cellular population is not too large.  Section \ref{section5} presents various numerical studies that support our theoretical findings.  Finally, we discuss our results and propose some open problems for future studies in Section \ref{section6}.

\vspace*{-10pt}
\section{Linearized stability analysis of homogeneous steady state}\label{section2}
In the mathematical analysis of pattern formation in reaction--diffusion systems, the principle of exchange of stability (\cite{CR4,Sat,Si} e.g.) is often employed to determine when bifurcation occurs for the family of evolution equations.  Generally speaking, this principle states that when a spatially homogeneous solution of a system becomes unstable as a parameter crosses a threshold value, it may admit stable spatially inhomogeneous solutions.  In particular, if the homogeneous solution loses stability through a pair of complex conjugate eigenvalues crossing the imaginary axis, one may expect, under suitable but reasonably technical conditions, that there exist time--periodic solutions to the evolution equations.  Moreover, this principle usually gives a qualitative relationship between the shape of bifurcating curve (such as its turning direction) of solutions and their stabilities.

In this paper, we are interested in studying time--periodic solutions to (\ref{11}), in contrast to the stable steady states investigated in \cite{WZYH}.  For the simplicity of our calculations and without much loss of generality, we shall confine our attention to system (\ref{11}) over one--dimensional interval $(0,L)$
\begin{equation}\label{21}
\left\{
\begin{array}{ll}
u_t=(d_1u_x-\chi u w_x)_x+\mu_1(1-u-a_1v)u,&x \in (0,L),t>0, \\
v_t=(d_2v_x-\xi v w_x)_x+\mu_2(1-a_2u-v)v,&x \in (0,L),t>0, \\
w_t=w_{xx}-\lambda w +u+v,&x \in (0,L),t>0, \\
u_x(x,t)=v_x(x,t)=w_x(x,t)=0,&x=0,L,t>0,\\
u(x,0)=u_0(x),v(x,0)=v_0(x),w(x,0)=w_0(x),&x\in(0,L).
\end{array}
\right.
\end{equation}
One of our primary goals is to explore how cellular kinetics effect the formation of spatially inhomogeneous positive solutions to (\ref{21}).  For this purpose, we adopt the principle of exchange of stability in the context of Hopf bifurcation, i.e., the bifurcation of a family of time--periodic solutions from $(\bar u, \bar v,\bar w)$.
To begin with, we carry out the linearized stability of $(\bar u,\bar v, \bar w)$ to investigate the spatial--temporal of dynamics to (\ref{21}) around this homogeneous solution.  Some of our stability results have been obtained in \cite{WZYH} and more details are needed for our purpose in this paper, therefore we include them here for the completeness and consistency of our arguments.

Linearizing (\ref{21}) by setting $U$, $V$ and $W$
\[u=\bar u+\epsilon U, v=\bar v+\epsilon V, w=\bar w+\epsilon W,\]
with $0<\epsilon \ll 1$ and substituting these perturbations into (\ref{21}), we obtain
\begin{equation}\label{22}
\left\{
\begin{array}{ll}
U_t\approx(d_1U'-\chi\bar u W')'-\mu_1\bar uU-\mu_1a_1\bar uV,&x \in (0,L),t>0,      \\
V_t\approx(d_2V'-\xi\bar v W')'-\mu_2a_2\bar vU-\mu_2\bar vV,&x \in (0,L),t>0,      \\
W_t\approx W''+U+V-\lambda W ,&x \in (0,L),t>0,      \\
U'(x)=V'(x)=W'(x)=0,&x=0,L,t>0.
\end{array}
\right.
\end{equation}
Now, we look for solutions of (\ref{22}) in the form $(U,V,W)=(C_1,C_2,C_3)e^{\sigma t+i\textbf{k}x}$, where $\textbf{k}$ is the wave mode vector and $\sigma$ is the growth rate of the perturbations respectively; here $\vert\textbf{k}\vert^2=(\frac{k\pi}{L})^2$ thanks to the $L^2$ eigen--expansions, and $C_i(i=1,2,3)$ are constants to be determined.    Substituting them into the linearized system above gives us the following problem
\[
\left(
\sigma I+\mathcal D\vert \textbf{k}\vert^2+\mathcal A_0
\right)
\left(
  \begin{array}{c}
C_1\\
C_2\\
C_3
\end{array}
\right)
=
\left(
  \begin{array}{c}
0\\
0\\
0
\end{array}
\right),
\]
where the matrices $\mathcal A_0$ and $\mathcal D$ are
\[
\mathcal A_0=\left(
  \begin{array}{ccc}
 -\mu_1\bar u & -\mu_1a_1\bar u & 0\\
    -\mu_2a_2\bar v & -\mu_2\bar v & 0\\
     1& 1 & -\lambda
\end{array}
\right),
\mathcal D=\left(
  \begin{array}{ccc}
d_1  & 0 & -\chi\bar u \\
0& d_2 & -\xi\bar v \\
0& 0 & 1
\end{array}
\right).
\]
Or equivalently $\sigma$ is an eigenvalue of the following stability matrix associated with (\ref{22})
\begin{equation}\label{23}
\mathcal A_k=
\left(
  \begin{array}{ccc}
    -d_1(\frac{k\pi}{L})^2-\mu_1\bar u & -\mu_1a_1\bar u & \chi\bar u(\frac{k\pi}{L})^2\\
    -\mu_2a_2\bar v & -d_2(\frac{k\pi}{L})^2-\mu_2\bar v & \xi\bar v(\frac{k\pi}{L})^2\\
     1& 1 & -(\frac{k\pi}{L})^2-\lambda
\end{array}
\right), k\in\mathbb N^+.
\end{equation}
By the standard principle of linearized stability (Theorem 5.2 in \cite{Si} or \cite{Sat} e.g.), $(\bar u,\bar v,\bar w)$ is asymptotically stable with respect to (\ref{21}) if and only if the real parts of all eigenvalues to matrix (\ref{23}) are negative.  The characteristic polynomial of (\ref{23}) is
\begin{equation}\label{24}
\sigma^3+\alpha_2(k)\sigma^2+\alpha_1(\chi,k)\sigma+\alpha_0(\chi,k)=0,
\end{equation}
where
\begin{equation*}
\alpha_2(k)= (d_1+d_2+1) \Big(\frac{k\pi}{L}\Big)^2+\mu_1\bar u+\mu_2\bar v+\lambda>0,
\end{equation*}
\begin{align*}
\alpha_1(\chi,k)=&\Big(\Big(\frac{k\pi}{L}\Big)^2+\lambda \Big) \Big(\big(d_1+d_2\big)\Big(\frac{k\pi}{L}\Big)^2
+\mu_1\bar u+\mu_2\bar v \Big)-a_1a_2\mu_1\mu_2\bar u\bar v \nonumber \\
&-(\chi\bar u+\xi\bar v)\Big(\frac{k\pi}{L}\Big)^2+\Big(d_1\Big(\frac{k\pi}{L}\Big)^2+\mu_1\bar u\Big) \Big(d_2\Big(\frac{k\pi}{L}\Big)^2+\mu_2 \bar v\Big),
\end{align*}
and
\begin{align*}
\alpha_0(\chi,k)&=\!-\chi\bar u\Big(\frac{k\pi}{L}\Big)^2\!\Big(d_2\Big(\frac{k\pi}{L}\Big)^2
\!\!+\!(1\!-\!a_2)\mu_2\bar v \Big)\!-\!\xi\bar v\Big(\frac{k\pi}{L}\Big)^2 \!\Big(d_1\Big(\frac{k\pi}{L}\Big)^2\!+\!(1\!-\!a_1)\mu_1\bar u \Big) \nonumber\\
&-\!a_1a_2\mu_1\mu_2\bar u\bar v\Big(\!\Big(\frac{k\pi}{L}\Big)^2\!\!\! +\!\lambda \Big)\!+\!\Big(d_1\Big(\frac{k\pi}{L}\Big)^2\!\!+\!\mu_1\bar u\Big)\Big(\!d_2\Big(\frac{k\pi}{L}\Big)^2\!\!+\!\mu_2\bar v\Big)\!\Big(\Big(\!\frac{k\pi}{L}\Big)^2\!\! +\!\lambda \Big).
\end{align*}
According to the Routh--Hurwitz conditions or Corollary 2.2 in \cite{LSW}, the real parts of all eigenvalues to (\ref{24}) are negative (hence $(\bar u,\bar v,\bar w)$ is locally stable) if and only if
\[\alpha_0(\chi,k)>0, \alpha_1(\chi,k)>0, \text{~and~}\alpha_1(\chi,k)\alpha_2(k)-\alpha_0(\chi,k)>0,\]
 for all $k\in \mathbb N^+$, while there exist some eigenvalues with a nonnegative real part if one of the conditions above fails for some $k\in \mathbb N^+$.  Moreover, since $\alpha_2(k)>0$, we will always have $\alpha_1(\chi,k)>0$ whenever $\alpha_0(\chi,k)>0$ and $\alpha_1(\chi,k)\alpha_2(k)-\alpha_0(\chi,k)>0$.  Therefore, the stability criterion above implies that $(\bar u,\bar v,\bar w)$ is unstable if there exists $k\in \mathbb N^+$ such that either $\alpha_0(\chi,k)<0$ or $\alpha_1(\chi,k)\alpha_2(k)-\alpha_0(\chi,k)<0$.  The following results are proved in \cite{WZYH}.
\begin{proposition}\label{proposition21}
Suppose that $0\leq a_1,a_2<1$ and all the rest parameters in (\ref{21}) are positive. Then the positive constant solution $(\bar u,\bar v,\bar w)$ of (\ref{21}) is unstable if $\chi\geq\chi_0=\min_{k\in \mathbb{N}^+}\{\chi^S_k,\chi^H_k\}$ and it is locally asymptotically stable if $\chi<\chi_0=\min_{k\in \mathbb{N}^+}\{\chi^S_k,\chi^H_k\}$,
where
\begin{align}\label{25}
\chi^S_k=&\frac{\big( \big(d_1(\frac{k\pi}{L})^2+\mu_1\bar u\big) \big(d_2(\frac{k\pi}{L})^2+\mu_2\bar v\big)-a_1a_2\mu_1\mu_2\bar u\bar v  \big) \big((\frac{k\pi}{L})^2+\lambda \big)}{d_2\bar u(\frac{k\pi}{L})^4+(1-a_2)\mu_2\bar u\bar v(\frac{k\pi}{L})^2}\nonumber \\
&-\frac{\xi \big(d_1\bar v(\frac{k\pi}{L})^4+(1-a_1)\mu_1\bar u\bar v(\frac{k\pi}{L})^2\big)}{d_2\bar u(\frac{k\pi}{L})^4+(1-a_2)\mu_2\bar u\bar v(\frac{k\pi}{L})^2},
\end{align}
and
\begin{equation}\label{26}
\chi^H_k=\frac{A_1{A_2}^2+{A_1}^2A_2+A_2A_3-\xi B_2}{B_1},
\end{equation}
with
\[A_1=(\frac{k\pi}{L})^2+\lambda,~A_2=(d_1+d_2)(\frac{k\pi}{L})^2+\mu_1\bar u+\mu_2\bar v,\]
\[A_3=\Big(d_1(\frac{k\pi}{L})^2+\mu_1\bar u\Big) \Big(d_2(\frac{k\pi}{L})^2+\mu_2\bar v \Big)-a_1a_2\mu_1\mu_2\bar u\bar v,\]
and
\[B_1=(d_1+1)\bar u(\frac{k\pi}{L})^4+(\lambda+\mu_1\bar u+a_2\mu_2\bar v)\bar u(\frac{k\pi}{L})^2,\]
\[B_2=(d_2+1)\bar v(\frac{k\pi}{L})^4+(\lambda+a_1\mu_1\bar u+\mu_2\bar v)\bar v(\frac{k\pi}{L})^2.\]
\end{proposition}
We want to point out that Proposition \ref{proposition21} holds for multi-dimensional bounded domains in $\mathbb{R}^N$ with $(\frac{k\pi}{L})^2$ being replaced by the $k$--th Neumann eigenvalue of $-\Delta$.

According to Proposition \ref{proposition21}, the spatially homogeneous steady state $(\bar u,\bar v,\bar w)$ loses its stability at $\chi_0=\min_{k\in \mathbb N^+}\{\chi^S_k,\chi^H_k\}$.  It is natural to expect that as $\chi$ surpasses the threshold value $\chi_0$, this homogeneous solution is driven unstable by spatially inhomogeneous solutions from the viewpoint of the principle of exchange of stability.  We use the indices $S$ and $H$ on the shoulder of $\chi_k$ to indicate that the stability is lost through steady state and Hopf bifurcation respectively as $\chi$ crosses $\chi^S_k$ and $\chi^H_k$.  One of the goals of this paper is to establish time--periodic patterns to (\ref{21}), in contrast to the stable steady states obtained in \cite{WZYH}.  Indeed, the authors in \cite{WZYH} carried out rigorous steady state bifurcation analysis on (\ref{21}) which shows that if $\chi_0=\min_{k\in\mathbb N^+}\chi^S_k$, the stability of $(\bar u,\bar v,\bar w)$ is lost to stable spatially inhomogeneous steady state of (\ref{21}); moreover, weakly nonlinear stability analysis near the bifurcating steady states is also performed which provides a wave mode selection mechanism.  On the other hand, numerical simulations in \cite{WZYH} suggest that (\ref{21}) admits stable time--periodic solutions when $\chi$ is around $\chi_0=\min_{k\in\mathbb N^+}{\chi^H_k}$.

Our approach is based on the Hopf bifurcation theorem for which \cite{A,CR4} and \cite{JD,Sat1} are good references.  For example, according to Theorem 1 in \cite{A} or Theorem 1.11 in \cite{CR4}, one of the necessary conditions for $\chi_k^H$ to be a bifurcation value of (\ref{21}) is that the stability matrix (\ref{23}) with $\chi=\chi_k^H$ has purely imaginary eigenvalues.  We first claim that $\chi_k^H\neq\chi_k^S$ if Hopf bifurcation occurs at $\chi=\chi_k^H$,  $\forall k\in\mathbb N^+$.  If not and we assume that $\chi_k^H=\chi_k^S$, then $\alpha_0(\chi,k)=\alpha_1(\chi,k)=0$ and (\ref{23}) has three eigenvalues $\sigma_1(k)=-\alpha_2(k)<0$, $\sigma_{2,3}(k)=0$, under which Hopf bifurcation does not occur.  To apply the Hopf bifurcation theorem in \cite{A,CR4}, we also need that, if $k\neq j$, matrices (\ref{23}) with $\chi=\chi_j^H$ and $\chi=\chi_k^H$ have different purely imaginary eigenvalues, which implies that $\chi_k^H\neq \chi_j^H$ if $k\neq j$.  Therefore we shall assume the following conditions in the rest of our analysis.
\begin{equation}\label{27}
\chi_k^H\neq \chi_k^S, \forall k\in \mathbb N^+, \text{~and~}\chi_k^H \neq \chi_j^H, \forall k\neq j \in\mathbb N^+.
\end{equation}
From straightforward calculations, we have that $\alpha_1(\chi,k)\alpha_2(k)=\alpha_0(\chi,k)$ if and only if $\chi=\chi^H_k$, and $\alpha_0(\chi,k)=0$ if and only if $\chi=\chi^S_k, \forall k\in \mathbb N^+$.  Therefore, if $\chi=\chi^S_{k}$, (\ref{24}) becomes $\sigma^3+\alpha_2(k)\sigma^2+\alpha_1(\chi_k^S,k)\sigma=0$, and if $\chi=\chi^H_{k}$, (\ref{24}) becomes $\sigma^3+\alpha_2(k)\sigma^2+\alpha_1(\chi_k^H,k)\sigma+\alpha_1(\chi_k^H,k)\alpha_2(k)=0$.  The following results are immediate from straightforward calculations.
\begin{proposition}\label{proposition22}
If $\chi=\chi^S_k$, i.e., $\alpha_0(\chi^S_k,k)=0$, (\ref{23}) has three eigenvalues given by $\tilde \sigma_1(k)=0$ and $\tilde \sigma_{2,3}(k)=\frac{-\alpha_2(k)\pm\sqrt{\alpha^2_2(k)-4\alpha_1(\chi^S_k,k)}}{2}$; if $\chi=\chi^H_k$, $\alpha_0(\chi^H_k,k)=\alpha_1(\chi^H_k,k)\alpha_2(k)$ and (\ref{23}) has three eigenvalues $\hat \sigma_1(k)=-\alpha_2(k)<0$ and $\hat \sigma_{2,3}(k)=\pm\sqrt{-\alpha_1(\chi^H_k,k)}$.
\end{proposition}
Hopf bifurcation occurs for (\ref{21}) in $(\bar u,\bar v,\bar w)$ only if (\ref{23}) has purely imaginary eigenvalues, and according to Proposition \ref{proposition22}, this is possible only when $\chi=\chi_k^H$ and $\alpha_1(\chi_k^H,k)>0$ since $\alpha_2(k)>0$, $\forall k\in\mathbb N^+$.  To determine when $\alpha_1(\chi_k^H,k)>0$, we denote $\bar \chi_k^1$ as the unique root of $\alpha_1(\chi,k)=0$ which is explicitly given by
\begin{align}\label{28}
\bar \chi_k^1=&\frac{\big((\frac{k\pi}{L})^2+\lambda\big)\big((d_1+d_2)(\frac{k\pi}{L})^2
+\mu_1\bar u+\mu_2\bar v \big)
+\big(d_1(\frac{k\pi}{L})^2+\mu_1\bar u\big) \big(d_2(\frac{k\pi}{L})^2+\mu_2 \bar v\big)}{\bar u(\frac{k\pi}{L})^2}\nonumber\\
&-\frac{\xi\bar v(\frac{k\pi}{L})^2+a_1a_2\mu_1\mu_2\bar u\bar v}{\bar u(\frac{k\pi}{L})^2}.
\end{align}
Now we have the following results.
\begin{lemma}\label{lemma21}
Let $\bar \chi_k^1$ be given by (\ref{28}).  Then for each $k\in\mathbb N^+$, we have that either (i) $\chi_k^H<\bar \chi_k^1<\chi_k^S$ or (ii) $\chi_k^S<\bar \chi_k^1<\chi_k^H$ occurs; moreover, if (i) occurs we have that $\alpha_1(\chi_k^H,k)>0>\alpha_1(\chi_k^S,k)$, and if (ii) occurs we have that $\alpha_1(\chi_k^S,k)>0>\alpha_1(\chi_k^H,k)$.
\end{lemma}
According to Lemma \ref{lemma21} and our discussions above, Hopf bifurcation may occur at $(\bar u,\bar v,\bar w,\chi_k^H)$ only when $\chi_k^H<\chi_k^S$.

\begin{remark}\label{remark21}
In general, it is very difficult to determine exactly when case \emph{(i)} or case \emph{(ii)} occurs in terms of system parameters.  However, if the interval length $L$ is sufficiently small, $\chi^S_k \approx \frac{d_1}{\bar {u}}(\frac{k\pi}{L})^2$ and $\chi^H_k\approx  \frac{(d_1+d_2)^2+(d_1+d_2)d_1d_2+d_1+d_2}{(d_1+1)\bar {u}}(\frac{k\pi}{L})^2$, therefore we always have that $\chi_k^S<\chi_k^H$, $\forall k\in \mathbb N^+$.  This fact indicates that (\ref{23}) has no purely imaginary eigenvalues when $L$ is sufficiently small.
\end{remark}


\section{Spatially inhomogeneous periodic patterns}\label{section3}
In this section, we prove the existence of time--periodic spatial patterns of (\ref{21}).  To be precise, we will show that, under proper assumptions on system parameters, the constant equilibrium $(\bar u,\bar v,\bar w)$ loses its stability through Hopf bifurcation as $\chi$ surpass $\chi_0=\min_{k\in \mathbb N^+}\{\chi_k^H,\chi_k^S\}$.  According to our analysis in Section \ref{section2}, the stability matrix (\ref{23}) has a paired purely imaginary eigenvalues if and only if $\chi=\chi^H_k<\chi_k^S$ and there does not exist a time--periodic solutions to (\ref{21}) that bifurcates from $(\bar u,\bar v,\bar w)$ if $\chi_k^H>\chi_k^S$.  Therefore, we assume that $\chi_k^H<\chi_k^S$ in the sequel in order to perform bifurcation analysis of (\ref{21}) at $\chi_k^H$.

\subsection{Hopf bifurcation}
For any $\chi>0$, we denote the eigenvalues of (\ref{23}) by $\sigma_1(\chi,k)$, $\sigma_2(\chi,k)$ and $\sigma_3(\chi,k)$.  When $\chi=\chi_k^H<\chi_k^S$, Proposition \ref{proposition22} and Lemma \ref{lemma21} tell us that eigenvalues of the stability matrix (\ref{23}) are: $\sigma_1(\chi_k^H,k)=-\alpha_2(k)$ and $\sigma_{2,3}(\chi_k^H,k)=\pm i\sqrt {\alpha_1(\chi_k^H,k)}$, which are purely imaginary.  Therefore, when $\chi$ is around $\chi_k^H$, (\ref{23}) has eigenvalues $\sigma_1(\chi)$ which is a real number around $-\alpha_2(k)$, and $\sigma_{2,3}(\chi,k)=\eta(\chi)\pm i\zeta(\chi)$, where $\eta(\chi)$ and $\zeta(\chi)$ are real analytical functions of $\chi$ satisfying $\eta(\chi_k^H)=0$ and $\zeta(\chi_k^H)=\sqrt{\alpha_1(\chi_k^H,k)}>0$.  In order to apply the bifurcation theory from \cite{A} or \cite{CR4} at point $\chi^{H}_{k}$, we need to verify the eigenvalue crossing condition or the so--called transversality condition.

Let us introduce the following Sobolev space
\[\mathcal X=\{u\in H^2(0,L) \vert u'(0)=u'(L)=0\}.\]
Our main result on the existence of nontrivial periodic orbits of (\ref{21}) states as follows.
\begin{theorem}\label{theorem31}
Assume that the parameters $d_i$, $\mu_i$, $a_i$, $i=1,2$, $\lambda$ and $\xi$ are positive and (\ref{27}) holds.  For each $k\in\mathbb N^+$, suppose that $0<\chi^{H}_k< \chi^{S}_{k}$ and $\chi^{H}_{j}\neq \chi^{H}_{k}$, $\forall j \neq k$, $j\in\mathbb N^+$.  Then there exist $\delta >0$ and a unique one--parameter family of nontrivial periodic orbits $\rho_k(s)=\Big(\textbf{u}_k(s, x, t), T_k(s), \chi_k(s) \Big)$: $s\in (-\delta,\delta) \rightarrow C^3(\mathbb{R},\mathcal X^3)\times \mathbb R^+\times\mathbb{R}$ satisfying $\Big( \textbf{u}_k(0, x, t),T_k(0),\chi_k(0) \Big)=\Big((\bar u,\bar v,\bar w), \frac{2\pi}{\zeta_0} ,\chi^{H}_{k} \Big)$
and
\begin{align}\label{31}
 \textbf{u}_k(s,x,t)=(\bar u,\bar v,\bar w)+s\Big(V^{+}_k e^{i\zeta_{0}t}+V^{-}_k e^{-i\zeta_{0}t}\Big) \cos \frac{k\pi x}{L}+o(s)
\end{align}
such that $(\textbf{u}_k(s),\chi_k(s))$ is a nontrivial   solution of (\ref{21}) and $\textbf{u}_k(s)$ is periodic with period
\begin{align}\label{32}
T_k(s)\approx \frac{2\pi}{\zeta_0}, \zeta_0=\sqrt{\alpha_1(\chi^{H}_k,k)}
\end{align}
and $\{V^{\pm}_k, \pm i\zeta_0\}$ are eigen--pairs of matrix $\mathcal A_k$; moreover $\rho_k(s_1)\neq \rho_k(s_2)$ for all $s_1 \neq s_2$, $\in(-\delta, \delta)$ and all nontrivial periodic solutions of (\ref{21}) around $(\bar u,\bar v,\bar w,\chi^H_k)$ must be on the orbit $\rho_k(s)$, $s\in(-\delta,\delta)$ in the sense that, if (\ref{21}) has a nontrivial periodic solution $\tilde {\textbf{u}}(x,t)$ with period $T$ for some $\chi\in\mathbb{R}$ around $\rho_k(s)$ such that $\vert \chi-\chi^{H}_{k}\vert <\epsilon$, $\vert T-2\pi/\zeta_0 \vert<\epsilon$ and $ \max_{t\in\mathbb{R}^+,x\in{\bar \Omega}}\vert \tilde U(x,t)-(\bar u,\bar v, \bar w) \vert<\epsilon$, where $\epsilon>0$ is a small constant, then there exist numbers $s\in(-\delta,\delta)$ and some $\theta\in [0,2\pi)$ such that $(T, \chi)=(T_k(s), \chi^H_k(s))$ and $\tilde {\textbf{u}}(x,t)=\textbf{u}_k(s,x,t+\theta)$.
\end{theorem}
\begin{proof}
Our proof is based on Theorem 1 from \cite{A} (or Theorem 1.11 from \cite{CR4}, Theorem 6.1 from \cite{LSW}).  Denote $\textbf{u}=(u,v,w)^T$.  We rewrite (\ref{21}) into the following abstract form
\[
\left\{
\begin{array}{ll}
\textbf{u}_t=(\mathcal D_0 \textbf{u}_x)_x+\textbf{f},&x\in(0,L), t>0,\\
\textbf{u}_x=0,&x=0,L, t>0,
\end{array}
\right.
\]
where
\begin{equation*}
\mathcal{D}_0 =\begin{pmatrix}
 d_1  &  0& -\chi u  \\
  0    & d_2& -\xi v\\
  0       & 0    &  1
  \end{pmatrix},
\textbf{f}=\begin{pmatrix}
\mu_1(1-u-a_1v)u\\
\mu_2(1-a_2u-v)v\\
-\lambda w+u+v
  \end{pmatrix},
\end{equation*}
then we see that system (\ref{21}) is normally parabolic since all the eigenvalues of $\mathcal{D}_0$ are positive.

Linearizing (\ref{21}) about $(\bar u,\bar v,\bar w)$ gives rise to the eigenvalue problem (\ref{22}), whose eigenvalues are those of the stability matrix $\mathcal A_k$ in (\ref{23}).  According to Proposition \ref{proposition22} and Lemma \ref{lemma21}, $\alpha_1(\chi_k^H)>0$ if $\chi=\chi_k^H<\chi_k^S$, and then matrix $\mathcal A_k$ in (\ref{23}) has a paired purely imaginary eigenvalues $\pm i \zeta_0=\pm i\sqrt{\alpha_1(\chi_k^H,k)}$.  Since $\chi_k^H \neq \chi_j^H$ for any $j\neq k$, matrix $\mathcal A_k$ has no eigenvalues of the form $iN\zeta_0$ for $N\in\mathbb N^+\backslash\{\pm1\}$; moreover 0 cannot be an eigenvalue for $\mathcal A_k$ with $\chi=\chi^H_k$ since $\chi_k^H<\chi_k^S$ according to Lemma \ref{lemma21}.

Let $\sigma_1(\chi)$, $\sigma_{2,3}(\chi)=\eta(\chi)\pm i\zeta(\chi)$ be the unique eigenvalue of (\ref{23}) in a neighbourhood of $\chi=\chi^H_k$, where we have skipped the index $k$ in each eigenvalue without confusing our reader.  It is easy to know that $\sigma_1$, $\eta$ and $\zeta$ are real analytical functions of $\chi$ with $\eta(\chi^{H}_{k})=0$ and $\zeta(\chi^{H}_{k})=\zeta_0>0$.  Theorem \ref{theorem31} follows from Theorem 1 of \cite{A} once we can prove the following \emph{transversality condition}
\[\frac{\partial \eta(\chi) }{\partial \chi}\Big\vert_{\chi=\chi^H_k}\neq0.\]
In particular, we will show that $\eta'(\chi^{H}_{k})>0$, where the prime $'$ here and in the sequel means the derivative taken with respect to $\chi$.  Substituting the eigenvalues $\sigma_1(\chi)$ and $\eta(\chi)\pm i\zeta(\chi)$ into (\ref{24}) and equating the real and imaginary parts there, we have that
\begin{align}\label{33}
-\alpha_2(\chi)&=2\eta(\chi)+\sigma_1(\chi), \\
\alpha_1(\chi,k)&=\eta^2(\chi)+\zeta^2(\chi)+2\eta(\chi)\sigma_1(\chi), \\
-\alpha_0(\chi,k)&=(\eta^2(\chi)+\zeta^2(\chi))\sigma_1(\chi).
\end{align}
Differentiating the equations above with respect to $\chi$, since $\alpha_2$ is independent of $\chi$, we obtain that
\begin{equation}\label{36}
2\eta'(\chi)+\sigma'_1(\chi)=0,
\end{equation}
and
\begin{align}\label{37}
&\left(
\begin{array}{ccc}
    \eta(\chi)-\sigma_1(\chi)  &\!\! 2\zeta(\chi)\\
    \eta^2(\chi)+\zeta^2(\chi)-\eta(\chi)\sigma_1(\chi) &\!\! 2\zeta(\chi)\sigma_1(\chi)
\end{array}
\right)
\!\!\left(
\begin{array}{ccc}\sigma'_1(\chi)\\
\zeta'(\chi)
\end{array}
\!\!\right) \nonumber\\
=& \left(
\begin{array}{ccc}-\bar u (\frac{k\pi}{L})^2\\
\bar u(\frac{k\pi }{L})^2  \big(d_2 (\frac{k\pi }{L})^2+\!(1\!-\!a_2)\mu_2\bar v \big)
\end{array}
\right).
\end{align}
Since $\eta(\chi^{H}_{k})=0$ and $\sigma_1(\chi^{H}_{k})=-\alpha_2(k)$, solving (\ref{37}) with $\chi=\chi_k^H$ gives us that
\[\sigma_1'(\chi^H_k)=-\frac{\bar u(\frac{k\pi }{L})^2 \big( (d_1+1)(\frac{k\pi }{L})^2+\mu_1\bar u+\mu_2 a_2\bar v+\lambda \big)}{\zeta_0^2+\alpha_2^2(k)}<0,\]
which implies $\eta'(\chi^H_k)=-\frac{1}{2}\sigma_1'(\chi^H_k)>0$ in light of (\ref{36}).  This verifies all the neces-\break sary conditions required in Theorem 1 of \cite{A}, from which our Theorem follows.
\end{proof}
Theorem \ref{theorem31} establishes time--periodic spatial patterns to (\ref{21}) bifurcating from $(\bar u,\bar v,\bar w)$.  Moreover it determines the exact bifurcation point $\chi_k^H$ and gives the explicit expression of the oscillation patterns, which admit spatial profile of the eigenfunction $\cos \frac{k\pi x}{L}$.  The arguments and results in Theorem \ref{theorem31} carry over to multi--dimensional bounded domain $\Omega$.

We have to point out that the condition $\chi_k^H<\chi_k^S$ is necessary for the occurrence of Hopf bifurcation at $(\bar u,\bar v,\bar w,\chi_k^H)$.  Considering the complexity of both terms, it is very hard to determine or evaluate when $\chi_k^H<\chi_k^S$, however according to Remark \ref{remark21}, for each $k\in\mathbb N^+$, if the interval length $L$ is sufficiently small, we must have that $\chi_k^S<\chi_k^H$.  This indicates that Hopf bifurcation does not occur in this situation hence (\ref{21}) does not have time--periodic solutions bifurcating from $(\bar u,\bar v,\bar w)$ when the interval length is sufficiently small.  Under this condition, it is showed in \cite{WZYH} that the stability of the homogeneous solution is lost through steady state bifurcation at the first bifurcation branch which has stable stationary solutions of (\ref{21}) with wave mode $\cos \frac{\pi x}{L}$.  See Theorem 3.2 in \cite{WZYH}.

It is worth mentioning that when $\mu_1=\mu_2=0$, we will show in Section \ref{section4} that $\chi_k^S<\chi_k^H$ for all $k\in\mathbb N^+$, regardless of the interval length.  This implies that $\chi_k^H$ is no longer a Hopf bifurcation point hence there does not exist any time--periodic solutions to (\ref{21}) that bifurcate from $(\bar u,\bar v, \bar w, \chi_k^H)$.

\subsection{Stability of time--periodic bifurcating solutions}
We proceed to analyze stability of the time--periodic bifurcating solutions on the bifurcation curves $\rho_k(s)$, $s\in(-\delta,\delta)$, obtained in Theorem \ref{theorem31}.  By stability here we mean the formal linearized stability of a periodic solution relative to disturbances from $\rho_k(s)$.  Assume that all the conditions in Theorem \ref{theorem31} are satisfied.  Suppose that $\chi^H_{k_0}=\min_{k\in\mathbb N^+}\chi^H_k<\chi_k^S$, $\forall k\in\mathbb N^+$, then our stability results show that $\rho_{k}(s)$, $s\in(-\delta,\delta)$, is asymptotically stable only if $k=k_0$, and $\rho_k(s)$, $s\in(-\delta,\delta)$, is always unstable for any $k\neq k_0$.  Certainly this is a necessary condition for stability.  Moreover, a rigorous mathematical treatment of stability away from the small--amplitude periodic solution is so far nonexistent.

Hopf's pioneering work in 1942 established the basic properties of time--periodic solutions to ODEs, such as existence and uniqueness, symmetry properties and stability, etc.  Since then, a considerable amount of work has been done in studying the stability of time--periodic solutions to Navier--Stokes equations \cite{Io,Iu,Sat} or abstract evolutions equations \cite{CR4,Jo,MM,Sat1}.  The stability of time--periodic solutions refers to the behavior of the Floquet multiplier or Floquet exponent \cite{HKW,Henry,Io,Iu} and we refer the reader to \cite{Bellman, Henry} or \cite{CR4} for reviews on the Floquet theory.

Denote $\textbf{u}_k(s,t)=(u_k(s,x,t), v_k(s,x,t), w_k(s,x,t))$ and let $(\textbf{u}_k(s,t),\chi_k(s))$ be the periodic solutions on the branch $\rho_k(s)$ obtained in Theorem \ref{theorem31}.  Rewrite (\ref{21}) into the following abstract form
\[\frac{d\textbf{u}_k}{dt}=\mathcal G(\textbf{u}_k,\chi_k(s)),\]
where
\[
\mathcal G(\textbf{u}_k,\chi_k(s))=\begin{pmatrix}
(d_1u_x-\chi_k(s) uw_x)_x+\mu_1(1-u-a_1v)u\\
(d_2v_x-\xi vw_x)_x+\mu_2(1-a_2u-v)v\\
 w_{xx}-\lambda w+u+v
\end{pmatrix}
\]
and we skip the index $k$ in $(u,v,w)$ without confusing our reader.  Differentiating the abstract system against $t$, writing $\dot {\textbf{u}}=\frac{d\textbf{u}}{dt}$, we have that
\[\frac{d \dot {\textbf{u}}_k}{dt}=\mathcal G_u(\textbf{u}_k, \chi_k(s))\dot {\textbf{u}}_k,\]
then we observe that $0$ is a Floquet exponent and 1 is a Floquet multiplier for $\textbf{u}_k$.

Linearize the periodic solution around the bifurcation branch $\rho_k(s)$ by substituting the perturbed solution $\textbf{u}_k+\textbf{w}e^{-\kappa t}$, where $\textbf{w}$ is a sufficiently small $T$--periodic function and $\kappa=\kappa(s)$ is a continuous function of $s$, then we have that
\begin{equation}\label{38}
\frac{d\textbf{w}(s,t)}{dt}=\mathcal G_u(\textbf{u}_k, \chi_k(s))\textbf{w}(s,t)+\kappa(s) \textbf{w}(s,t) ,
\end{equation}
where $\mathcal G_u$ is the Fr\'echet derivative with respect to $\textbf{u}$.  Then stability of the bifurcating solutions in the neighborhood of the branch point $\chi_k^H$ can be determined by computing the eigenvalues of this reduced equation.  At $s=0$ (\ref{38}) is associated with the eigenvalue problem
\begin{equation}\label{39}
\mathcal G_0(k)\textbf{w}=\kappa(0) \textbf{w},
\end{equation}
where
\[\mathcal G_0(k)=\mathcal G_u((\bar u,\bar v,\bar w), \chi_k^H)=
\left(
  \begin{array}{ccc}
    d_1\frac{d^2}{dx^2}-\mu_1\bar u & -\mu_1a_1\bar u & -\chi_k^H\bar u\frac{d^2}{dx^2}\\
    -\mu_2a_2\bar v & d_2\frac{d^2}{dx^2}-\mu_2\bar v & -\xi\bar v\frac{d^2}{dx^2}\\
     1& 1 & \frac{d^2}{dx^2}-\lambda
\end{array}
\right).
\]
It is easy to see that the spectrum of $\mathcal G_0$ is infinitely dimensional; in particular, we want to point out that $\mathcal G_0$ corresponds to the stability matrix of $(\bar u,\bar v,\bar w)$ given in (\ref{23}).
\begin{equation}\label{310}
\mathcal A_j(\chi_k^H)=
\left(
  \begin{array}{ccc}
    -d_1(\frac{j\pi}{L})^2-\mu_1\bar u & -\mu_1a_1\bar u & \chi_k^H \bar u(\frac{j\pi}{L})^2\\
    -\mu_2a_2\bar v & -d_2(\frac{j\pi}{L})^2-\mu_2\bar v & \xi\bar v(\frac{j\pi}{L})^2\\
     1& 1 & -(\frac{j\pi}{L})^2-\lambda
\end{array}
\right), j\in\mathbb N^+.
\end{equation}
Suppose that $\min_{k\in\mathbb N^+} \{\chi_k^H,\chi_k^S\}=\chi^H_{k_0}$ for some $k_0\in\mathbb N^+$.  We first show that $\rho_{k}(s)$ around $\chi^H_k$ is unstable for any $k\neq k_0$.

Indeed, denote the eigenvalues to $\mathcal A_{k}(\chi_k^H)$ by $\sigma_1(\chi_k^H, k)$, $\sigma_2(\chi_k^H, k)$ and $\sigma_3(\chi_k^H, k)$, then we have that the real part of one of these eigenvalues must be positive, i.e., Re$(\sigma_2(\chi_k^H, k))>0$, since $(\bar u,\bar v,\bar w)$ is unstable if $\chi>\chi_0$ according to Proposition \ref{proposition21}. Therefore for any positive integer $k\neq k_0$, we have that $\mathcal G_0(k)$ must have an eigenvalue with positive real part, hence $\kappa(0)<0$ if $k\neq k_0$.  By the standard perturbation theory for an eigenvalue of finite multiplicity (see \cite{Henry} or \cite{Kato} e.g.), $\kappa(s)<0$ for $s$ being small if $k\neq k_0$, therefore all the bifurcation branches $\rho_k(s)$ around $(\bar u,\bar v,\bar w)$ are unstable if $k\neq k_0$.  This implies that if a periodic bifurcating solution is stable, it must be on the $k_0$--th branch at which $\chi_k^H$ achieves its minimum over $\mathbb N^+$ (i.e., it is on the left--most branch), while all the later branches are always unstable.

\begin{figure}[h!]
        \centering
        \begin{subfigure}[b]{0.45\textwidth}
                \includegraphics[width=\textwidth]{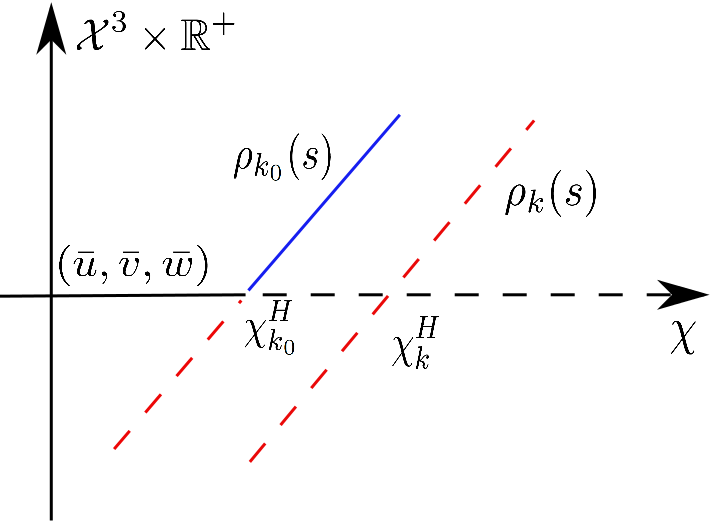}
                \caption*{Supercritical bifurcation $\chi^{H'}_{k}(0)>0$.}
        \end{subfigure}\hspace{0.15in}
        \begin{subfigure}[b]{0.45\textwidth}
        \includegraphics[width=\textwidth]{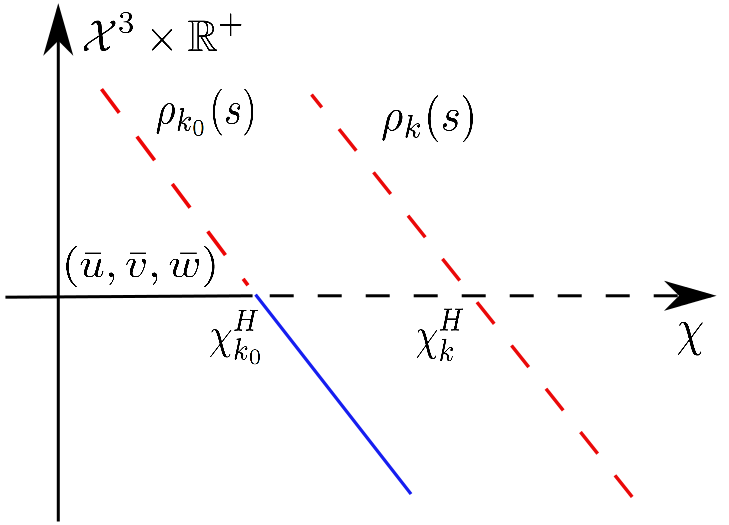}
                \caption*{Subcritical bifurcation $\chi^{H'}_{k}(0)<0$.}
        \end{subfigure}
 \caption{Bifurcation diagrams of $\rho_k(s)$ around $(\bar u,\bar v,\bar w)$.  The stable bifurcation curve is plotted in solid lines  and the unstable bifurcation curve is plotted in imaginary line.  The branch $\rho_{k}(s)$ around $(\bar u,\bar v,\bar w,\chi_{k})$ is always unstable if $k\neq k_0$, while the turning direction of $\rho_{k_0}(s)$ determines its stability.}\label{figure1}
\end{figure}

We proceed to discuss stability of branch $\rho_{k_0}(s)$ around $(\bar u,\bar v,\bar w,\chi_{k_0}^H)$.  According to Lemma 2.10 in \cite{CR4}, the eigenvalue $\kappa(s)$ is a continuous real function of $s$ which is uniquely defined near $s=0$.  For $\chi$ being around $\chi_{k_0}^H$, the eigenvalues of $\mathcal A_{k_0}$ are $\sigma_1(\chi)$, $\sigma_{2,3}(\chi)=\eta(\chi)\pm i\zeta(\chi)$.  By Theorem 2.13 in \cite{CR4}, $\kappa(s)$ and $s\chi'_{k_0}(s)$ have the same zeros in small neighbourhood of $s=0$ in which $\kappa(s)$ and $-\eta'(\chi_{k_0}^H)s\chi'_{k_0}(s)$ have the same sign (if they are not zero), and
\[\vert \kappa(s)+\eta'(\chi_{k_0}^H\big)s\chi'_{k_0}(s)\vert\leq \vert s \chi'_{k_0}(s) \vert o(1), \text{~as~}s\rightarrow 0.\]
By Theorem 8.2.3 in \cite{Henry}, the bifurcating periodic solutions are stable if $\kappa(s)>0$, and they are unstable if $\kappa(s)<0$.  Moreover, since we already showed in the proof of Theorem \ref{theorem31} that $\eta'(\chi_{k_0}^H)>0$, $\kappa(s)$ has the same sign as $s\chi'_{k_0}(s)$, therefore if $\chi''_{k_0}(0)\neq 0$, the branching solutions are stable if they appear supercritically and unstable if they appear subcritically.  See Theorem 3 in the survey paper \cite{Sat1} of D. H. Sattinger.  The stability of bifurcation branches around $(\bar u,\bar v,\bar w)$ is schematically presented in Figure \ref{figure1}.
In order to evaluate $\chi'_{k_0}(0)$, one can follow the calculations using the factorization theorem in \cite{JN} or the method of integral averaging in \cite{CM}, or the normal form method and center manifold theorem from \cite{HKW} of Hassard \emph{et al.}.  In each method, we need to perform a perturbation analysis in the neighbourhood of the critical bifurcation value $\chi^H_{k_0}$, by substituting $\chi^H_k(s)$ and the periodic solution $\textbf{u}(s,x,t)$ as Taylor series of $s$ into (\ref{21}), then we equate the $s$--terms to find algebraic equations of $\chi'_{k_0}(0)$ which determines the direction of the Hopf bifurcation if $\chi''_{k_0}(0)\neq 0$.  The calculations are routine but extremely complicated, therefore we skip them here.

\section{System without cellular growth}\label{section4}
In this section, we study the positive solutions to (\ref{21}) with $\mu_1=\mu_2=0$, i.e., the following system
\begin{equation}\label{41}
\left\{
\begin{array}{ll}
u_t=(d_1u_x-\chi u w_x)_x ,&x \in (0,L),t>0     \\
v_t=(d_2v_x-\xi v w_x)_x ,&x \in (0,L),t>0      \\
w_t=w_{xx}-\lambda w +u+v,&x \in (0,L),t>0      \\
u_x(x,t)=v_x(x,t)=w_x(x,t)=0,&x=0,L,t>0,        \\
u(x,0)=u_0(x),v(x,0)=v_0(x),w(x,0)=w_0(x),&x\in(0,L).
\end{array}
\right.
\end{equation}
One of our main results in this section shows that (\ref{41}) and its multi--dimensional counterpart have no positive time--periodic solutions and it indicates that cellular growth is responsible for the formation of time--periodic spatial patterns in system (\ref{11}) and (\ref{21}).  According to our results in Section \ref{section3} and \cite{WZYH}, $(\bar u,\bar v,\bar w)$ loses its stability to steady state bifurcating solutions when $\chi_0=\min_{k\in\mathbb N^+} \chi^S_k$ and to Hopf bifurcating solutions when $\chi_0=\min_{k\in\mathbb N^+} \chi^H_k$.  We shall show that Hopf bifurcation does not occur for (\ref{21}) at $(\bar u,\bar v,\bar w)$ when $\mu_1=\mu_2=0$ from the viewpoint of linearized stability analysis and then we proceed to investigate the effect of cellular kinetics on the dynamics of (\ref{41}).  In particular, we shall show that the kinetics are necessary for the formation of periodic patterns of the two competing species chemotaxis system.

\subsection{Linearized stability of $(\bar u,\bar v,\bar w)$}
Similar as in Section \ref{section2}, our starting point is the linearized stability analysis of the constant solution $(\bar u,\bar v,\bar w)$.  To this end, we perform some elementary calculations to show that Hopf bifurcation never takes place at $(\bar u,\bar v,\bar w,\chi_k^H)$ if $\mu_1=\mu_2=0$, thanks to which the stability matrix (\ref{23}) becomes
\begin{equation}\label{42}
\mathcal A_k=
\left(
  \begin{array}{ccc}
    -d_1(\frac{k\pi}{L})^2 & 0 & \chi\bar u(\frac{k\pi}{L})^2\\
    0 & -d_2(\frac{k\pi}{L})^2 & \xi\bar v(\frac{k\pi}{L})^2\\
     1& 1 & -(\frac{k\pi}{L})^2-\lambda
\end{array}
\right).
\end{equation}
By the same arguments that lead to Proposition \ref{proposition21}, we have that $(\bar u,\bar v,\bar w)$ is unstable with respect to (\ref{41}) if $\chi\geq\chi_0=\min_{k\in \mathbb{N}^+} \{\chi^S_k, \chi^H_k\}$ and it is locally asymptotically stable if $\chi<\chi_0$ where $\chi_k^S$ in (\ref{25}) and $\chi_k^H$ in (\ref{26}) become
\[\chi^S_k=\frac{d_1d_2\big((\frac{k\pi}{L})^2+\lambda\big)-\xi d_1\bar v}{d_2 \bar u},\]
\[\chi_k^H=\frac{A_1{A^{2}_2+A^{2}_1A_2+A_2A_3-\xi B_2}}{B_1}\]
with
\[A_1=(\frac{k\pi}{L})^2+\lambda,~~A_2=(d_1+d_2)(\frac{k\pi}{L})^2,~~A_3=d_1d_2(\frac{k\pi}{L})^4,\]
and
\[B_1=(d_1+1)\bar u (\frac{k\pi}{L})^4+\lambda\bar u(\frac{k\pi}{L})^2,~~B_2=(d_2+1)\bar v (\frac{k\pi}{L})^4+\lambda\bar v(\frac{k\pi}{L})^2.\]
Similar as above, since we consider positive chemotaxis (chemical being chemo--attractive to the cells), we must have that $\chi_0>0$ hence $\chi_k^S>0$, which implies
\begin{align}\label{43}
0<\xi<\frac{d_2\big((\frac{k\pi}{L})^2+\lambda\big)}{\bar v}.
\end{align}
To see that there does not exist time--periodic positive solutions to (\ref{41}) bifurcating from $(\bar u, \bar v,\bar w)$,  we shall prove that $\chi_k^H$ is not a bifurcation value by showing that $\chi^S_k<\chi^H_k$ for all $k\in \mathbb{N}^+$.

To this end, we denote $\chi_k^S-\chi_k^H=\frac{J_1+J_2}{J_3}$, where
\[J_1=\Big(d_1d_2\Big(\Big(\frac{k\pi}{L}\Big)^2+\lambda\Big)-\xi d_1\bar v\Big)\Big(\big(d_1+1\big)\Big(\frac{k\pi}{L}\Big)^2+\lambda\Big)+\xi d_2\bar v\Big(\big(d_2+1\big)\Big(\frac{k\pi}{L}\Big)^2+\lambda\Big),\]
\begin{align*}
    &J_2=-\big(d_1+d_2\big)^2d_2\Big(\frac{k\pi}{L}\Big)^2\Big(\Big(\frac{k\pi}{L}\Big)^2+\lambda\Big)-\big(d_1+d_2\big)d_2\Big(\Big(\frac{k\pi}{L}\Big)^2 +\lambda\Big)^2\\&\quad\quad -\big(d_1+d_2\big)d_1d^{2}_2\Big(\frac{k\pi}{L}\Big)^4,
\end{align*}
and
\[J_3=\big(d_1+1\big)d_2\bar u\Big(\frac{k\pi}{L}\Big)^2+\lambda d_2\bar u.\]
Since $J_3>0$, we just need to determine the sign of $J_1+J_2$ in order to find that of $\chi_k^S-\chi_k^H$.  Simple calculations give rise to
\begin{align*}
&J_1+J_2\nonumber\\
=&\Big(d_1d_2\Big(\Big(\frac{k\pi}{L}\Big)^2+\lambda\Big)-\xi d_1\bar v\Big)\Big(\big(d_1+1\big)\Big(\frac{k\pi}{L}\Big)^2+\lambda\Big)+\xi d_2\bar v\Big(\big(d_2+1\big)\Big(\frac{k\pi}{L}\Big)^2+\lambda\Big)\nonumber\\
&-\big(d_1+d_2\big)^2d_2\Big(\frac{k\pi}{L}\Big)^2\Big(\Big(\frac{k\pi}{L}\Big)^2+\lambda\Big)-\big(d_1+d_2\big)d_2\Big(\Big(\frac{k\pi}{L}\Big)^2+\lambda\Big)^2\\
& -\big(d_1+d_2\big)d_1d^{2}_2\Big(\frac{k\pi}{L}\Big)^4\nonumber\\
=&d_1d_2\Big(\Big(\frac{k\pi}{L}\Big)^2+\lambda\Big)\Big(\big(d_1+1\big)\Big(\frac{k\pi}{L}\Big)^2+\lambda\Big)-\xi\bar v\Big(d_1\big((d_1+1)\Big(\frac{k\pi}{L}\Big)^2+\lambda\big)\\
&-d_2\big((d_2+1)\Big(\frac{k\pi}{L}\Big)^2+\lambda\big)\Big)-\big(d_1+d_2\big)^2d_2\Big(\frac{k\pi}{L}\Big)^2\Big(\Big(\frac{k\pi}{L}\Big)^2+\lambda\Big)\nonumber\\
&-(d_1+d_2)d_2\Big(\Big(\frac{k\pi}{L}\Big)^2+\lambda\Big)^2-(d_1+d_2)d_1d^{2}_2\Big(\frac{k\pi}{L}\Big)^4\nonumber,
\end{align*}
and we divide our discussions into the following two cases:

\vspace*{4pt}
\noindent\textbf{Case 1}.  If $d_1\geq d_2$, then since $\xi>0$ from (\ref{43}), we have
\begin{align*}
&J_1+J_2\\
<&d_1d_2\Big(\Big(\frac{k\pi}{L}\Big)^2+\lambda\Big)\Big(\big(d_1+1\big)\Big(\frac{k\pi}{L}\Big)^2+\lambda\Big)-(d_1+d_2)^2d_2\Big(\frac{k\pi}{L}\Big)^2\Big(\Big(\frac{k\pi}{L}\Big)^2+\lambda\Big)\nonumber\\
&-(d_1+d_2)d_2\Big(\Big(\frac{k\pi}{L}\Big)^2+\lambda\Big)^2-(d_1+d_2)d_1d^{2}_2\Big(\frac{k\pi}{L}\Big)^4\nonumber\\
=&d_2\Big(\Big(\frac{k\pi}{L}\Big)^2+\lambda\Big)\Big(\big(d_1(d_1+1)-(d_1+d_2)^2-(d_1+d_2)\big)\Big(\frac{k\pi}{L}\Big)^2\\&+\big(d_1-(d_1+d_2)\big)\lambda\Big)
-(d_1+d_2)d_1d^{2}_2\Big(\frac{k\pi}{L}\Big)^4\nonumber\\
=&-d^{2}_2\Big(\Big(\frac{k\pi}{L}\Big)^2+\lambda\Big)\Big((2d_1+d_2+1)\Big(\frac{k\pi}{L}\Big)^2+\lambda\Big)-(d_1+d_2)d_1d^{2}_2\Big(\frac{k\pi}{L}\Big)^4\nonumber\\
<&0.
\end{align*}

\vspace*{4pt}
\noindent\textbf{Case 2}.  If $d_1<d_2$, we apply the fact $\xi<\frac{d_2((\frac{k\pi}{L})^2+\lambda)}{\bar v}$ in (\ref{43}) to estimate
\begin{align*}
&J_1+J_2\\
<&d_1d_2\Big(\Big(\frac{k\pi}{L}\Big)^2+\lambda\Big)\Big(\big(d_1+1\big)\Big(\frac{k\pi}{L}\Big)^2+\lambda\Big)+d_2\Big(\Big(\frac{k\pi}{L}\Big)^2+\lambda\Big)\\
&\Big(d_2\big((d_2+1)\Big(\frac{k\pi}{L}\Big)^2+\lambda\big)
-d_1\big((d_1+1)\Big(\frac{k\pi}{L}\Big)^2+\lambda\big)\Big)\\
&-\big(d_1+d_2\big)^2d_2\Big(\frac{k\pi}{L}\Big)^2\Big(\Big(\frac{k\pi}{L}\Big)^2+\lambda\Big)-(d_1+d_2)d_2\nonumber\\
&\cdot\Big(\Big(\frac{k\pi}{L}\Big)^2+\lambda\Big)^2-(d_1+d_2)d_1d^{2}_2\Big(\frac{k\pi}{L}\Big)^4,\nonumber   \\
<&d^{2}_2\Big(\Big(\frac{k\pi}{L}\Big)^2+\lambda\Big)\Big(\big(d_2+1\big)\Big(\frac{k\pi}{L}\Big)^2+\lambda\Big)-(d_1+d_2)^2d_2\Big(\frac{k\pi}{L}\Big)^2\Big(\Big(\frac{k\pi}{L}\Big)^2+\lambda\Big)\nonumber\\
&-(d_1+d_2)d_2\Big(\Big(\frac{k\pi}{L}\Big)^2+\lambda\Big)^2-(d_1+d_2)d_1d^{2}_2\Big(\frac{k\pi}{L}\Big)^4\nonumber\\
=&d_2\Big(\Big(\frac{k\pi}{L}\Big)^2+\lambda\Big)\Big(  \big( d_2(d_2+1)-(d_1+d_2)^2-(d_1+d_2)\big)\Big(\frac{k\pi}{L}\Big)^2\\
&+\big(d_2-(d_1+d_2)\big)\lambda \Big)
-(d_1+d_2)d_1d^{2}_2\Big(\frac{k\pi}{L}\Big)^4,\nonumber   \\
=&-d_1d_2\Big(\Big(\frac{k\pi}{L}\Big)^2+\lambda\Big)\Big((d_1+2d_2+1)\Big(\frac{k\pi}{L}\Big)^2+\lambda\Big)-(d_1+d_2)d_1d^{2}_2\Big(\frac{k\pi}{L}\Big)^4\nonumber\\
<&0,
\end{align*}
therefore we have that $J_1+J_2<0$ in both cases, hence $\chi_k^S<\chi_k^H$ for each $k\in\mathbb N^+$ as claimed.  According to Proposition \ref{proposition22} and Lemma \ref{lemma21}, matrix (\ref{42}) does not have purely imaginary eigenvalues for any $\chi$ and therefore (\ref{41}) has no Hopf bifurcating solutions from $(\bar u,\bar v,\bar w)$.
\begin{remark}\label{remark41}
If $\Omega$ is a bounded domain in $\mathbb R^N$, $N\geq2$, then the constant solution is unstable if $\chi>\chi_0$, with $(\frac{k\pi}{L})^2$ being replaced by the $k$--th Neumann eigenvalue of $-\Delta$.  Without loss of generality, we assume that $(\bar u, \bar w,\bar w)$ is the same as given by (\ref{12}).  If not, then thanks to the conservation of cellular populations, we must have that
\[(\bar u, \bar v)=\frac{1}{\vert \Omega\vert}\Big(\int_\Omega u_0, \int_\Omega v_0\Big),\]
and our calculations above still hold true under the new notations.
\end{remark}

According to our discussions above, $(\bar u,\bar v,\bar w)$ loses stability to steady state bifurcating solutions as $\chi$ surpasses $\chi_0=\min_{k\in\mathbb N^+}\{\chi_k^H,\chi_k^S\}$.  However, since $\chi_k^S<\chi_k^H$, we have from Proposition \ref{proposition22} and Lemma \ref{lemma21} that the stability matrix (\ref{42}) has no purely imaginary eigenvalue, therefore Hopf bifurcation can not occur for system (\ref{41}), which does not admit stable time--periodic patterns bifurcating from $(\bar u,\bar v,\bar w)$.  Moreover, according to the results in \cite{WZYH}, we know that the stability of $(\bar u,\bar v,\bar w)$ is lost through steady state bifurcation to spatially inhomogeneous patterns of (\ref{41}), which has a spatial profile $\cos\frac{k_0 \pi x}{L}$.

\subsection{Lyapunov functional}
The linearized stability analysis of $(\bar u,\bar v,\bar w)$ suggests that system (\ref{41}) has no time--periodic patterns that bifurcate from this constant solution and it does not rule out the existence of time--periodic patterns of (\ref{41}) since there may exist time oscillating solutions other than those from Hopf bifurcation.  However, we shall prove that the latter case is indeed impossible by showing the existence of time--monotone Lyapunov functional to (\ref{41}).  Our results also hold for (\ref{41}) over multi--dimensional domains hence we consider the following fully parabolic system
\begin{equation}\label{44}
\left\{
\begin{array}{ll}
u_t=\nabla \cdot (d_1 \nabla u-\chi u \nabla w),&x \in \Omega,t>0, \\
v_t=\nabla \cdot(d_2\nabla v-\xi v \nabla w),&x \in \Omega,t>0, \\
w_t=\Delta w-\lambda w +u+v,&x \in \Omega,t>0,\\
\frac{\partial u}{\partial \textbf{n}}=\frac{\partial v}{\partial \textbf{n}}=\frac{\partial w}{\partial \textbf{n}}=0,&x\in\partial \Omega,t>0,\\
u(x,0)=u_0(x),v(x,0)=v_0(x),w(x,0)=w_0(x), &x\in \Omega,
\end{array}
\right.
\end{equation}
where $\Omega\subset \mathbb R^N$, $N\geq1$, $\nabla $ is the gradient operator and $\Delta$ is the Laplace operator.  The system parameters are the same as in (\ref{21}).

We shall show that (\ref{44}) has a time--monotone Lyapunov functional, therefore it admits no time--periodic patterns regardless of space dimension and system parameters as long as there is no cellular growth.  Thanks to this fact, another main contribution of this paper is the global existence and large--time behavior of positive solutions to (\ref{44}).  We begin with the verification that (\ref{44}) has a Lyapunov functional in the following form
\begin{equation}\label{45}
F(u,v,w)=\frac{1}{2}\int_{\Omega} \vert \nabla w \vert^2 + \frac{\lambda}{2}\int_{\Omega} w^2 +\frac{d_1}{\chi}\int_{\Omega}(u\ln u-u)+\frac{d_2}{\xi}\int_{\Omega} (v\ln v-v) -\int_{\Omega} (u+v)w
\end{equation}
which is non--increasing along the trajectories of (\ref{44}).  We have the following Lemma.
\begin{lemma}\label{lemma41}
Suppose that $(u,v,w)$ is a classical solution of (\ref{44}) in $\Omega \times (0,T)$, $T\in(0,\infty]$ and the initial data $u_0$ and $v_0$ are strictly positive on $\bar \Omega$.  Then $F(u,v,w)$ given in (\ref{45}) is monotone decreasing in time and it satisfies
\begin{equation}\label{46}
\begin{split}
&\int_0^t\int_{\Omega} w^2_t +\int_0^t\int_{\Omega} \Big(\frac{( d_1\nabla u-\chi u\nabla w)^2}{\chi u}+\frac{( d_2\nabla v-\xi v\nabla w)^2}{\xi v}\Big)+F(u,v,w)\\
=&F(u_0,v_0,w_0).
\end{split}
\end{equation}
\end{lemma}
\begin{proof}
According to the Maximum Principles (e.g. \cite{LSU}) and positivity of the initial data, both $u$ and $v$ are strictly positive on $\bar \Omega \times (0,T)$.  We have from straightforward calculations that
\begin{align}\label{47}
\frac{dF}{dt}=& -\int_{\Omega} w_t(w_t+\lambda w-u-v)+\lambda\int_{\Omega} w w_t +\frac{d_1}{\chi}\int_{\Omega} u_t\ln u +\frac{d_2}{\xi}\int_{\Omega} v_t\ln v \nonumber\\
&-\int_{\Omega} (u_t+v_t)w -\int_{\Omega} (u+v)w_t  \nonumber \\
=& -\int_{\Omega} w^2_t -\int_{\Omega} (u_t+v_t)w+\frac{d_1}{\chi}\int_{\Omega}u_t\ln u+ \frac{d_2}{\xi}\int_{\Omega} v_t\ln v.
\end{align}
To estimate (\ref{47}), we have from the PDEs and the divergence theorem that
\begin{align}\label{48}
-\int_{\Omega} u_t w &=-\int_{\Omega} \nabla \cdot (d_1\nabla u-\chi u\nabla w)w= d_1\int_{\Omega}\nabla u\cdot \nabla w-\chi\int_{\Omega}u\vert \nabla w \vert^2,
\end{align}
and
\begin{align}\label{49}
-\int_{\Omega} v_t w =d_2\int_{\Omega}\nabla v\cdot\nabla w-\xi \int_{\Omega}v\vert \nabla w \vert^2,
\end{align}
while the last two terms of (\ref{47}) becomes
\begin{equation}\label{410}
\frac{d_1}{\chi}\int_{\Omega} u_t\ln u=\frac{d_1}{\chi}\int_{\Omega}\nabla \cdot(d_1\nabla u-\chi u\nabla w)\ln u=-\frac{d_1^2}{\chi}\int_{\Omega}\frac{\vert\nabla u\vert^2}{u} +d_1\int_{\Omega} \nabla u\cdot \nabla w,
\end{equation}
and
\begin{align}\label{411}
\frac{d_2}{\xi}\int_{\Omega} v_t\ln v=-\frac{d_2^2}{\xi}\int_{\Omega}\frac{\vert\nabla v\vert^2}{v} +d_2\int_{\Omega}\nabla v\cdot \nabla w .
\end{align}
In light of (\ref{48})--(\ref{411}), (\ref{47}) leads us to
\begin{align}\label{412}
\frac{dF}{dt}=&-\int_{\Omega} w^2_t-\frac{d_1^2}{\chi}\int_{\Omega} \frac{\vert \nabla u \vert^2}{u} + 2d_1\int_{\Omega}\nabla u\cdot\nabla w-\chi\int_{\Omega} u\vert \nabla w \vert^2 \nonumber\\
&- \frac{d_2^2}{\xi}\int_{\Omega}\frac{ \vert \nabla v \vert^2}{v} + 2d_2\int_{\Omega}\nabla v\cdot\nabla w-\xi \int_{\Omega}v\vert \nabla w \vert^2 \nonumber\\
=&-\int_{\Omega} w^2_t-\int_{\Omega} \Big(\frac{d_1}{\sqrt{\chi}}\frac{\nabla u}{\sqrt{u}}-\sqrt{\chi}\sqrt{u}\nabla w\Big)^2-\int_{\Omega} \Big(\frac{d_2}{\sqrt{\xi}}\frac{\nabla v}{\sqrt{v}}-\sqrt{\xi}\sqrt{v}\nabla w\Big)^2\nonumber\\
=&-\int_{\Omega} w^2_t-\int_{\Omega} \frac{( d_1\nabla u-\chi u\nabla w)^2}{\chi u}-\int_{\Omega}  \frac{( d_2\nabla v-\xi v\nabla w)^2}{\xi v}\leq 0.
\end{align}
Therefore $F(u,v,w)$ is always non--increasing in $t$ and (\ref{46}) follows from (\ref{412}).
\end{proof}
The time--monotone Lyapunov functional (\ref{45}) indicates that (\ref{44}) can not have time--periodic solutions, in contrast to system (\ref{21}) which has oscillating solutions according to Theorem \ref{theorem31}.  This indicates that the formation of oscillating solutions to (\ref{11}) is driven by the appearance of cellular growth terms.

\subsection{Global existence and boundedness for $N=2$}
In \cite{BEG}, the authors investigated parabolic--parabolic--elliptic system of (\ref{44}) with $\lambda=0$ over the whole space $\mathbb R^N$, $N\geq2$.  Their results state that, in loose terms, if the initial data $u_0$ and $v_0$ concentrate at some points $x_i$, $i\in \mathbb N^+$, then the solutions to (\ref{44}) can blow up within finite time.  In \cite{EVC}, Espejo \emph{et al}. studied the parabolic--parabolic--elliptic system of (\ref{44}) with $\lambda=d_2=1$ over a unit disk in $\mathbb R^2$ under homogeneous Dirichlet boundary conditions.  They showed that if
\[\int_\Omega u_0<\frac{8\pi d_1}{\chi}, \int_\Omega v_0<\frac{8\pi}{\xi} \text{ and }\Big(\int_\Omega u_0+\int_\Omega v_0\Big)^2<\frac{8\pi d_1}{\chi}\int_\Omega u_0+\frac{8\pi}{\xi}\int_\Omega v_0,\]
then there exist global bounded classical positive solutions.  In \cite{CEV}, it is proved that if one of the inequalities above fails, then the solutions to a similar problem over $\mathbb R^2$ blow up.  For $\Omega=\mathbb R^2$, global existence and large--time behaviors for (\ref{44}) are investigated in \cite{ZL} provided that $\Vert (u_0,v_0,\nabla w_0) \Vert_{L^1(\mathbb R^2)}$ are sufficiently small, following the arguments on invariant sets of (\ref{44}) as in \cite{Winkler}.

This section is devoted to studying the global existence and boundedness of classical positive solutions to (\ref{44}) as well as their large--time behaviors.  Similar as for (\ref{11}), Amann's theories \cite{Am1,Am2} guarantee the local existence of (\ref{44}) since it is a \emph{normally parabolic triangle} system, while $L^1$--boundedness of $u$ and $v$
 still holds for (\ref{44}) due to the conservation of cellular populations.  By the same arguments for Theorem 2.5 in \cite{WZYH} we can prove the local existence and boundedness for (\ref{44}) over $(0,T_\text{max})$ for some $T_\text{max} \in (0,\infty]$.  We are mainly concerned with the global existence and boundedness of (\ref{44}) over $\Omega \subset \mathbb R^2$.  In particular, assuming
\begin{equation}\label{413}
\frac{\chi}{d_1}\int_{\Omega} u_0+\frac{\xi}{d_2}\int_{\Omega} v_0<4 \pi,
\end{equation}
we show that the positive classical solutions to (\ref{44}) exist globally and are uniformly bounded in time as follows.
\begin{theorem}\label{theorem42}
Let $\Omega\subset \mathbb R^2$ be a smooth and bounded domain.  Assume that the initial data $(u_0, v_0, w_0)\in W^{1,p}\times W^{1,p} \times C(\bar \Omega)$ for some $p>2$, and $u_0,v_0>0$ and $w_0\geq0$, $\not \equiv 0$ in $\Omega$.  Under condition (\ref{413}), system (\ref{44}) admits a unique classical global solution $(u(x,t), v(x,t), w(x,t))$ for all $(x,t)\in\Omega \times (0,\infty)$; moreover, $(u(x,t), v(x,t), w(x,t))$ is nonnegative in $\Omega \times (0,\infty)$ and $\Vert(u,v,w) \Vert_{L^\infty(\Omega)}$ is uniformly bounded for all $t\in (0,\infty)$.
\end{theorem}
In Figure \ref{figure2}, we plot the numerical simulations to illustrate the evolution of spatially--inhomogeneous time--periodic patterns of (\ref{21}).  The numerics there indicate the lack of a stable global attractor to the full system (\ref{21}), at least for the parameter set we choose.  Therefore, we are motivated to investigate large--time behavior of positive solutions to (\ref{44}) by establishing the existence time--monotone Lyapunov functional.  We show that the classical solutions to (\ref{44}) converge to its stationary states as time goes to infinity.  See Theorem \ref{theorem46}.  For example, the first subgraph in Figure \ref{figure4} plots the spatial--temporal dynamics of (\ref{44}) over $\Omega=(0,6)$, where the interior spike is an attractor of the system.

We now pass to present our proof of Theorem \ref{theorem42}.  Our main vehicle is the $L^2$--estimates proved in Lemma \ref{lemma45} and the application of standard Moser--Alikakos iteration \cite{A0}.  To derive the $L^2$ estimates, we shall estimate energy--type functionals $\Vert u\ln u\Vert_{L^1}$ and $\Vert v\ln v\Vert_{L^1}$ via a special version of the Moser--Trudinger inequality.  It is necessary to remark that the crucial use of the embedding inequalities only applies when $N\leq2$.  The following result is well known (e.g. \cite{CY}).
\begin{lemma}\label{lemma43}
Let $\Omega\subset \mathbb R^2$ be a smooth and bounded domain, then there exists a positive constant $C>0$ dependent on $\Omega$ such that for all $\phi\in H^1(\Omega)$
\begin{equation}\label{414}
\int_{\Omega} \exp \vert \phi\vert \le C \exp\Big\{\frac{1}{8\pi}
\Vert \nabla \phi\Vert^2_{L^2}+\frac{1}{\vert\Omega\vert} \Vert \phi\Vert_{L^1}\Big\}.
\end{equation}
\end{lemma}
Let $(u,v,w)$ be the classical positive solutions to (\ref{44}) over $\Omega \times (0,T_{\max})$, $T_{\max} \in (0,\infty]$, then by analogous arguments for Theorem 2 in \cite{BN} or Lemma 3.4 in \cite{NSY}, we can prove the following results.
\begin{lemma}\label{lemma44}
Assume that all the conditions in Theorem \ref{theorem42} are satisfied.  Then there exists a constant $C>0$ such that the solutions to (\ref{44}) have the property
\begin{equation}\label{415}
\int_{\Omega} (u+v)w <C \text{ and } \int_\Omega u \ln u+v\ln v<C, \forall t\in (0,T_{\max}).
\end{equation}
\end{lemma}
\begin{proof}
Since $x\rightarrow -\ln x$ is convex and $\int_{\Omega} \frac{u}{\int_{\Omega} u_0} =1$, $\int_{\Omega} \frac{v}{\int_{\Omega} v_0} =1$, by Jensen's inequality we have that for any $\delta>0$
\begin{align*}
-\ln \Big(\frac{1}{\int_{\Omega} u_0}\int_{\Omega} e^{\frac{(1+\delta)\chi}{d_1}w }\Big)&=-\ln \Big(\int_{\Omega} \frac{e^{\frac{(1+\delta)\chi}{d_1}w }}{u} \frac{u}{\int_{\Omega} u_0}\Big)\nonumber\\
&\le \int_{\Omega} \Big(-\ln \frac{e^{\frac{(1+\delta)\chi}{d_1}w }}{u}\Big) \frac{u}{\int_{\Omega} u_0}\nonumber\\
&=-\frac{1+\delta}{\int_{\Omega} u_0}\frac{\chi}{d_1}\int_{\Omega} uw+\frac{1}{\int_{\Omega} u_0} \int_{\Omega} u\ln u;
\end{align*}
multiplying this inequality by $\frac{d_1}{\chi}\int_\Omega u_0$ gives rise to
\begin{align}\label{416}
&\frac{d_1}{\chi}\int_{\Omega} u\ln u -\int_{\Omega} uw  \nonumber\\
\ge& \delta\int_{\Omega} uw+\frac{d_1}{\chi}\Big(\int_{\Omega} u_0\Big)\ln \Big(\int_{\Omega} u_0\Big)-\frac{d_1}{\chi}\Big(\int_{\Omega} u_0\Big) \ln\Big(\int_{\Omega} e^{\frac{(1+\delta)\chi}{d_1}w  }\Big).
\end{align}
Similarly, we can have from the $v$--equation that
\begin{align}\label{417}
&\frac{d_2}{\xi}\int_{\Omega} v\ln v -\int_{\Omega} vw  \nonumber\\
\ge& \delta\int_{\Omega} vw+\frac{d_2}{\xi}\Big(\int_{\Omega} v_0\Big)\ln \Big(\int_{\Omega} v_0\Big)-\frac{d_2}{\xi}\Big(\int_{\Omega} v_0\Big) \ln\Big(\int_{\Omega} e^{\frac{(1+\delta)\xi}{d_2}w}\Big).
\end{align}
On the other hand, we apply (\ref{414}) on (\ref{416}) and (\ref{417}) with $\phi=\frac{(1+\delta)\chi_i}{d_i}w$, $i=1,2$, respectively, where $\chi_1=\chi$ and $\chi_2=\xi$,  to have that
\begin{align}\label{418}
\ln \Big(\int_{\Omega} e^{\frac{(1+\delta)\chi_i}{d_i}w  }\Big)&\le \ln C +\frac{1}{8\pi}\int_{\Omega} \Big\vert \nabla \frac{(1+\delta)\chi_i}{d_i}w \Big\vert^2 +\frac{1}{\vert \Omega\vert}\int_\Omega \frac{(1+\delta)\chi_i}{d_i} w  \nonumber\\
&\le \frac{\chi_i^2(1+\delta)^2}{8\pi d_i^2} \int_{\Omega} \vert \nabla w \vert^2+C_i, i=1,2,
\end{align}
where we have applied the boundedness of $\Vert w \Vert_{L^1}$ in (\ref{418}).  In light of (\ref{416})--(\ref{418}), we have that
\begin{equation}\label{419}
F(u,v,w)\ge \frac{1}{2} \int_{\Omega} \vert \nabla w \vert^2 -\frac{(1+\delta)^2}{8\pi}\Big( \frac{\chi}{d_1}\int_{\Omega} u_0+\frac{\xi}{d_2}\int_{\Omega} v_0\Big)\int_{\Omega} \vert \nabla w \vert^2+\delta\int_{\Omega} (u+v)w-C_3,
\end{equation}
where $C_3$ is a positive constant that depends on $\Vert u \Vert_{L^1}+\Vert v \Vert_{L^1}$.  Choosing $\delta>0$ to be sufficiently small, we see from condition (\ref{413}) that
\[\frac{(1+\delta)^2}{8\pi}\Big(\frac{\chi}{d_1}\int_{\Omega} u_0+\frac{\xi}{d_2}\int_{\Omega} v_0\Big)\le \frac{1}{2},\]
which, together with (\ref{419}), implies that
\begin{equation}\label{420}
F(u,v,w)\ge \delta\int_{\Omega} (u+v)w-C_4.
\end{equation}
Now we can easily see that (\ref{420}) implies $\int_{\Omega} (u+v)w<C_5$.

On the other hand, we have from straightforward calculations
\begin{align}\label{421}
&\frac{d_1}{\chi}\int_{\Omega} u\ln u +\frac{d_2}{\xi}\int_{\Omega} v\ln v \nonumber\\
=&F(u,v,w)- \frac{1}{2}\int_{\Omega} \vert \nabla w \vert^2 - \frac{\lambda}{2}\int_{\Omega} w^2 +\frac{d_1}{\chi}\int_{\Omega}u  +\frac{d_2}{\xi}\int_{\Omega}v +\int_{\Omega} (u+v)w \nonumber\\
\le &F(u_0,v_0,w_0)+\frac{d_1}{\chi}\int_{\Omega}u  +\frac{d_2}{\xi}\int_{\Omega}v +\int_{\Omega} (u+v)w \le C_6,
\end{align}
therefore both $\int_\Omega u \ln u$ and $\int_\Omega v\ln v$ are uniformly bounded from above.  This completes the proof of Lemma \ref{lemma44}.
\end{proof}
The following result is an immediate consequence of (\ref{415}) with $p=3$ in Lemma 3.5 of \cite{NSY}.
\begin{corollary}\label{corollary1}
Assume the same conditions in Theorem \ref{theorem42}.  Let $u,v$ be the classical solutions to (\ref{44}), then for any $\epsilon>0$, there exists a positive constant $C(\epsilon)$ such that
\begin{equation}\label{422}
\Vert u\Vert^3_{L^3}\leq \epsilon \Vert \nabla u\Vert^2_{L^2}+C(\epsilon) \text{ and }\Vert v\Vert^3_{L^3}\leq \epsilon \Vert \nabla v\Vert^2_{L^2} +C(\epsilon).
\end{equation}
\end{corollary}

Next we provide the boundedness of $\Vert u(\cdot,t)\Vert_{L^2}+\Vert v(\cdot,t)\Vert_{L^2}$ for $t\in(0,\infty)$, which suffices to prove the global existence and boundedness of $(u,v,w)$ to (\ref{44}).
\begin{lemma}\label{lemma45}
Under the same conditions in Theorem \ref{theorem42}, there exists a positive constant $C$ such that
\begin{equation}\label{423}
\Vert u \Vert_{L^2(\Omega)}+\Vert v \Vert_{L^2(\Omega)}\le C, \forall t\in(0,T_{\max}).
\end{equation}
\end{lemma}
\begin{proof}In light of the PDEs in (\ref{44}), straightforward calculations involving integration by parts and Young's inequality lead us to
\begin{align}\label{424}
&\frac{1}{2}\frac{d}{dt} \int_{\Omega} u^2 =-d_1\int_{\Omega} \vert \nabla u\vert^2 +\chi\int_{\Omega} u\nabla u\cdot\nabla w =-d_1\int_{\Omega} \vert \nabla u\vert^2 +\frac{\chi}{2}\int_{\Omega} \nabla (u^2) \cdot\nabla w \nonumber\\
&=-d_1\int_{\Omega} \vert \nabla u\vert^2 -\frac{\chi}{2}\int_{\Omega} u^2\Delta w =-d_1 \int_{\Omega}\vert \nabla u\vert^2 -\frac{\chi}{2}\int_{\Omega} u^2(w_t+\lambda w-u-v)\nonumber\\
&\le -d_1 \int_{\Omega}\vert \nabla u\vert^2  -\frac{\chi}{2}\int_{\Omega}  u^2w_t +\frac{\chi}{2}\int_{\Omega}  u^3 + \frac{\chi}{2}\int_{\Omega}  u^2v \nonumber\\
&\le -d_1 \int_{\Omega}\vert  \nabla u\vert^2  -\frac{\chi}{2}\int_{\Omega}  u^2w_t +\frac{5\chi}{6}\int_{\Omega}  u^3 + \frac{\chi}{6}\int_{\Omega}  v^3 ,
\end{align}
and
\begin{equation}\label{425}
\frac{1}{2}\frac{d}{dt} \int_{\Omega} v^2 \le -d_2 \int_{\Omega}\vert  \nabla v\vert^2-\frac{\xi}{2}\int_{\Omega}  v^2w_t +\frac{5\xi}{6}\int_{\Omega}  v^3 + \frac{\xi}{6}\int_{\Omega}  u^3 .
\end{equation}
To estimate (\ref{424}) (similarly (\ref{425})), we have from the Gagliardo--Ladyzhenskaya--Nirenberg interpolation inequality (see \cite{LSU} e.g.) and Cauchy--Schwartz that, for any $u\in W^{1,4}(\Omega)$, there exist two positive constants $C_1$ and $C_2$ such that
\[\Vert u\Vert_{L^4}^2\le C_1\Vert \nabla u\Vert_{L^2}\Vert u\Vert_{L^2}+C_2 \Vert u\Vert^2_{L^2}.\]
Moreover we have from H\"older's that in (\ref{424})
\begin{align}\label{426}
&-\frac{\chi}{2}\int_{\Omega}  u^2w_t \le \frac{\chi}{2}  \Vert u^2\Vert_{L^2} \Vert w_t\Vert_{L^2}= \frac{\chi}{2}\Vert u\Vert_{L^4}^2 \Vert w_t\Vert_{L^2}\nonumber\\
\le& C_3\Vert \nabla u\Vert_{L^2}\Vert  u\Vert_{L^2}\Vert w_t\Vert_{L^2}+C_4\Vert u\Vert_{L^2}^2\Vert w_t\Vert_{L^2}\nonumber\\
\le& \frac{d_1}{4} \Vert \nabla u\Vert_{L^2}^2+ \frac{C^2_3}{d_1}\Vert  u\Vert_{L^2}^2\Vert w_t\Vert_{L^2}^2+C_4\Vert  u\Vert_{L^2}^2\Big(\frac{d_1}{4} +\frac{1}{d_1}\Vert w_t\Vert_{L^2}^2\Big)\nonumber\\
\le& \frac{d_1}{4} \Vert \nabla u\Vert_{L^2}^2+\frac{C_4 d_1}{4} \Vert  u\Vert_{L^2}^2+ \frac{C^2_3+C_4}{d_1} \Vert  u\Vert_{L^2}^2\Vert w_t\Vert_{L^2}^2,
\end{align}
and $C_3=\frac{\chi}{2}C_1$ and $C_4=\frac{\chi}{2}C_2$; similarly
\begin{align}\label{427}
-\frac{\xi}{2}\int_{\Omega}v^2w_t \le \frac{d_2}{4} \Vert \nabla v\Vert_{L^2}^2+\frac{C_5 d_2}{4} \Vert v\Vert_{L^2}^2+ \frac{C_5+C^2_6}{d_2} \Vert  v\Vert_{L^2}^2\Vert w_t\Vert_{L^2}^2,
\end{align}
where $C_5$ and $C_6$ are positive constants; moreover, we can have from (\ref{422}) that
\begin{equation}\label{428}
\frac{5\chi+\xi}{6}\int_\Omega u^3 \leq \frac{d_1}{4} \Vert \nabla u\Vert^2_{L^2} +C_7 \text{ and }  \frac{\chi+5\xi}{6}\int_\Omega v^3  \leq \frac{d_2}{4} \Vert \nabla v\Vert^2_{L^2} +C_8.
\end{equation}

Adding up (\ref{424})--(\ref{425}), using (\ref{426}), (\ref{427}) and (\ref{428}), we have
\begin{align}\label{429}
\frac{1}{2}\frac{d}{dt} \int_{\Omega} u^2 +\frac{1}{2}\frac{d}{dt} \int_{\Omega} v^2\le & -\frac{d_1}{2}\Vert  \nabla u\Vert_{L^2}^2-\frac{d_2}{2}\Vert  \nabla v\Vert_{L^2}^2 +C_9 \Vert  u\Vert_{L^2}^2+C_{10}\Vert v\Vert_{L^2}^2 \nonumber\\
&+ C_{11} \Vert  u\Vert_{L^2}^2\Vert w_t\Vert_{L^2}^2+ C_{12} \Vert  v\Vert_{L^2}^2\Vert w_t\Vert_{L^2}^2+C_7+C_8.
\end{align}
Moreover we have from Corollary 1 in \cite{CKWW} due to Gagliardo--Ladyzhenskaya-- \break Nirenberg inequality that for any $\epsilon>0$, there exist $C_{13}$ and $C_{14}$ such that
\[\Vert u\Vert_{L^2}^2\le \epsilon \Vert\nabla u\Vert_{L^2}^2+C_{13}(\epsilon) \text{ and } \Vert v\Vert_{L^2}^2\le \epsilon \Vert\nabla v\Vert_{L^2}^2+C_{14}(\epsilon),\]
where $C_{13}$ ($C_{14}$) is a positive constant that depends on $\epsilon$, $\Omega$ and $\Vert u_0\Vert_{L^1}$ ($\Vert v_0\Vert_{L^1}$).  Now we have from (\ref{429}) that
\begin{align}\label{430}
\frac{1}{2}\frac{d}{dt} \int_{\Omega} u^2 +\frac{1}{2}\frac{d}{dt} \int_{\Omega} v^2\le & \Big(C_9-\frac{d_1}{2\epsilon}\Big)\Vert u\Vert_{L^2}^2+\Big(C_{10}-\frac{d_2}{2\epsilon}\Big)\Vert v\Vert_{L^2}^2\nonumber\\
&+\Big(C_{11} \Vert  u\Vert_{L^2}^2+ C_{12} \Vert  v\Vert_{L^2}^2\Big)\Vert w_t\Vert_{L^2}^2+C_{15}(\epsilon).
\end{align}
Choosing
\[0<\epsilon<\frac{1}{2}\min\Big\{\frac{d_1}{2C_9}, \frac{d_2}{2C_{10}}\Big\},\]
and denoting
\[y(t)=\int_\Omega u^2(x,t)+v^2(x,t)dx,\]
we have from (\ref{430}) that
\begin{equation}\label{431}
y'(t)\le \Big(-\frac{d}{2\epsilon}+C_{16}\Vert w_t\Vert_{L^2}^2\Big)y(t)+C_{17}(\epsilon), \forall t\in(0,T_{\max}),
\end{equation}
where $d=\min\{d_1,d_2\}$, $C_{16}=\max\{2C_{11},2C_{12}\}$ is a positive constant independent of $\epsilon$ and $C_{17}$ is also a positive constant.

To derive the boundedness of $y(t)$, we recall from (\ref{46}) and (\ref{415}) that $\int_0^t \Vert w_t(\cdot, \break t)\Vert_{L^2}ds$ is bounded for all $t\in(0,\infty)$.  Solving (\ref{431}) gives rise to
\begin{align}\label{432}
y(t)\le& y_0\exp \Big\{\int_0^t \Big(-\frac{d}{2\epsilon}+C_{16}\Vert w_t\Vert_{L^2}^2 \Big)d\tau\Big\}\nonumber \\
 &+ C_{17}(\epsilon)\int_0^t \exp\Big\{\int_s^t \Big(-\frac{d}{2\epsilon}+C_{16}\Vert w_t\Vert_{L^2}^2 \Big) d\tau\Big\} ds,\nonumber \\
\le& C_{18} y_0 \exp \Big\{-\frac{d}{2\epsilon}t\Big\}+C_{19}(\epsilon),
\end{align}
where $y_0=\Vert u_0\Vert_{L^2}^2+\Vert v_0\Vert_{L^2}^2$ and $C_{18}$ is a positive constant that may depend on $\epsilon$.  This finishes the proof of Lemma \ref{lemma45}.
\end{proof}
\begin{proof}[Proof\nopunct] \emph{of Theorem} \ref{theorem42}.
The proof is exact the same as that of Theorem 2.5 in \cite{WZYH}, where Moser--Alikakos iteration and standard bootstrap arguments are applied, except that the global existence there is established for $N=1$, therefore we shall only sketch the main steps here.

First of all, since (\ref{44}) is a triangular system, its local existence follows from the classical results of Amann \cite{Am1,Am2} and the regularity of the solutions follows from standard parabolic regularity arguments.  By the same estimates for (2.7) in \cite{WZYH}, we can find a constant $C>0$ such that
\[
\begin{split}
\Vert w (\cdot,t) \Vert_{W^{1,q}(\Omega)} \leq& C\Big(1+ \int_0^t e^{-\nu(t-s)} (t-s)^{-\frac{1}{2}-\frac{2}{2}(\frac{1}{2}-\frac{1}{q})}\cdot\big(\Vert u(\cdot,s) \Vert_{L^2(\Omega)}
\\&\left. +\Vert v(\cdot,s) \Vert_{L^2(\Omega)}+\Vert w(\cdot,s) \Vert_{L^2(\Omega)}\big) ds\right.\Big),\\
\leq&  C\Big(1+ \sup_{s\in(0,t)} (\Vert (u,v,w)(\cdot,s) \Vert_{L^2(\Omega)}\Big),
\end{split}
\]
therefore $\Vert w (\cdot,t) \Vert_{W^{1,q}(\Omega)}$ is uniformly bounded for each $q\in(1,\infty)$ since $\Vert (u,v)(\cdot,\break t)\Vert_{L^2(\Omega)}$ are bounded.  Then we can again apply Gagliardo--Nirenberg interpolation to show the boundedness of $\Vert (u,v)(\cdot,t)\Vert_{L^3}$, which implies that $\Vert w (\cdot,t) \Vert_{W^{1,\infty}(\Omega)}$ is uniformly bounded.  Finally, after applying the standard Moser--Alikakos iteration, we obtain the uniform boundedness of $(u,v)$ in $L^\infty$.
\end{proof}

\subsection{Asymptotic behaviors and nonconstant positive steady states}
In this subsection, we will show that the bounded classical solutions to (\ref{44}) converge to the (probably nontrivial) steady state as $t\rightarrow \infty$.  Using the monotonicity of the Lyapunov functional $F$ in (\ref{45}), we can prove that the limit of $\mathcal \omega$--sets of (\ref{44}) is the stationary system.  The following theorem can be proved by the same arguments for Lemma 3.1 in \cite{Winkler} thanks to Lemma \ref{lemma45}.
\begin{theorem}\label{theorem46}
Let $\Omega$ be a bounded domain in $\mathbb R^2$.  Suppose that all conditions in Theorem \ref{theorem42} are satisfied.  Let $(u,v,w)$ be the global bounded classical solutions to (\ref{44}), then there exists $t_k\rightarrow \infty$ as $k\rightarrow \infty$ such that $(u(\cdot,t_k),v(\cdot,t_k),w(\cdot,t_k)) \rightarrow (u_\infty,v_\infty,w_\infty)$ in $C^2(\bar \Omega) \times C^2(\bar \Omega) \times C^2(\bar \Omega)$, where $(u_\infty,v_\infty,w_\infty)$ satisfies the following system
\begin{equation}\label{433}
\left\{
\begin{array}{ll}
\nabla \cdot (d_1 \nabla u_\infty-\chi u_\infty \nabla w_\infty)=0,&x \in \Omega \\
\nabla \cdot(d_2\nabla v_\infty-\xi v_\infty \nabla w_\infty)=0,&x \in \Omega, \\
\Delta w_\infty-\lambda w_\infty +u_\infty+v_\infty=0,&x \in \Omega,\\
\frac{\partial u_\infty}{\partial \textbf{n}}=\frac{\partial v_\infty}{\partial \textbf{n}}=\frac{\partial w_\infty}{\partial \textbf{n}}=0,&x\in\partial \Omega,\\
\int_\Omega u_\infty=\int_\Omega u_0, \int_\Omega v_\infty=\int_\Omega v_0.
\end{array}
\right.
\end{equation}
\end{theorem}
According to Remark \ref{remark41}, (\ref{433}) has no nonconstant stable steady state if $\chi$ is large.  It is interesting to study nonconstant positive solutions to (\ref{433}).  In particular, the steady states with concentrating properties such as boundary or interior spikes can be used to model the aggregation phenomenon for chemotactic cells.

\section{Numerical simulations}\label{section5}
In this section, we perform some numerical studies of stable and time--periodic spatially inhomogeneous solutions to system (\ref{21}).  To manifest the effect of cellular growth and other parameters on its spatial--temporal dynamics, we fix $a_1=a_2=0.5$ in all our simulations, thanks to which the IBVP has positive equilibrium $(\bar u,\bar v,\bar w)=(\frac{2}{3},\frac{2}{3},\frac{4}{3\lambda})$, while all initial data are selected to be small perturbations from the equilibrium.  We shall choose different sets of system parameters to study the initiation and development of spatial patterns to the system.

\begin{figure}[h!]
        \centering
\includegraphics[width=\textwidth,height=2.7in]{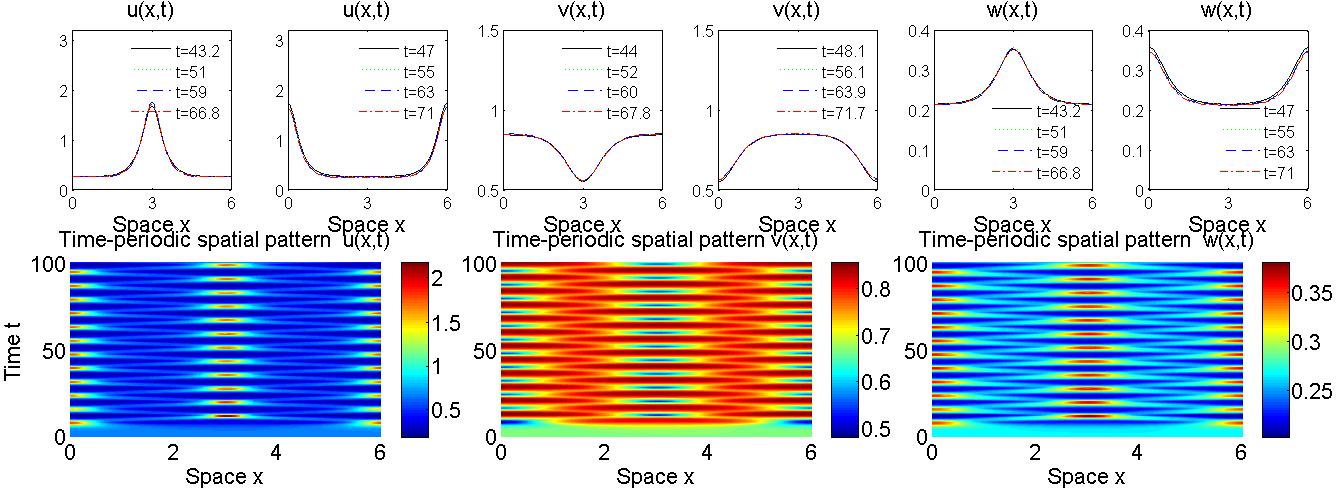}
  \caption{Initiation and development of time--periodic spatial patterns to (\ref{11}) over $(0,6)$ with initial data being small perturbations of $(\bar u,\bar v,\bar w)$.  System parameters are chosen to be $d_1=5$, $d_2=0.1$, $\mu_1=\mu_2=1$, $\lambda=5$, $\xi=0.1$ and $\chi=80$.  Our theoretical results indicate that the homogeneous equilibrium loses its stability at $\chi_0=\chi^H_{2}\approx 63.2$ through Hopf bifurcation to a stable time--periodic pattern which has spatial profile $\cos \frac{\pi x}{3}$ and period $T\approx 8$. Space and time grid sizes are $\Delta x=0.02$ and $\Delta t=0.05$.  The numerical simulations are in good agreement with our theoretical findings.}\label{figure2}
\end{figure}

\begin{figure}[h!]
        \centering
\includegraphics[width=\textwidth,height=3in]{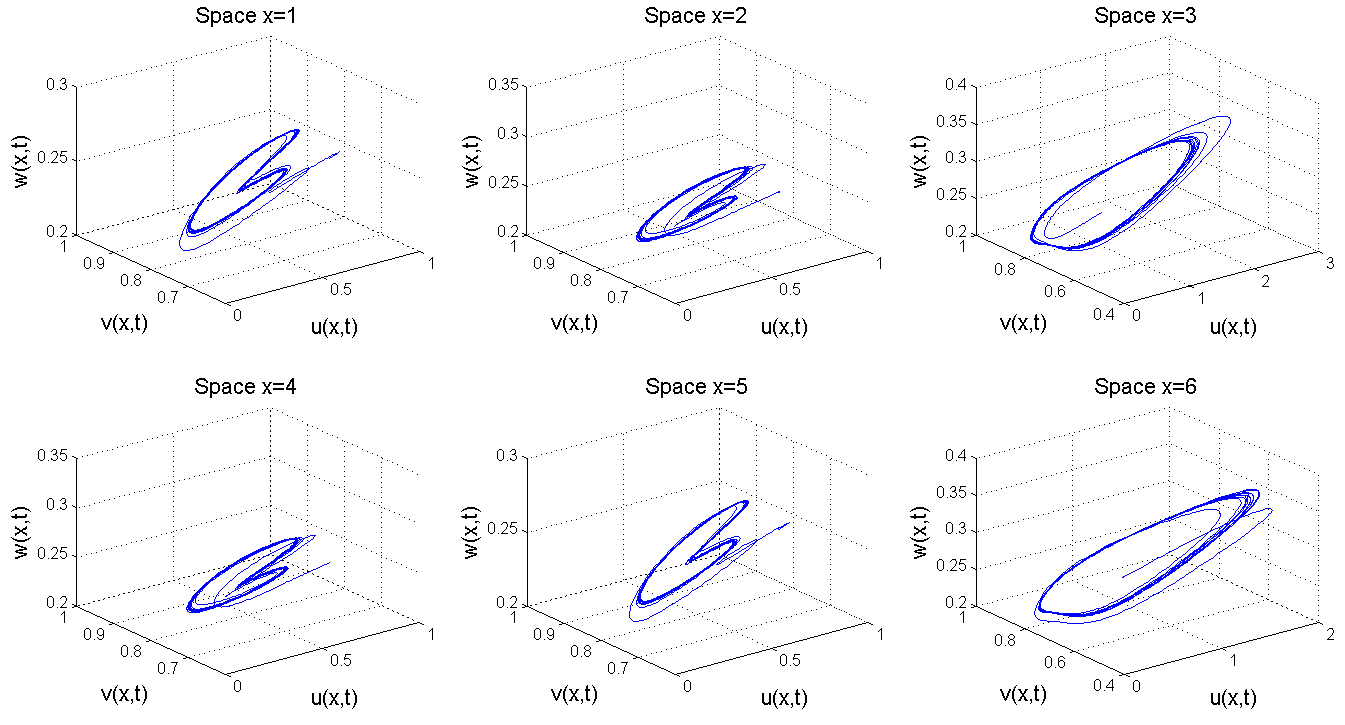}
  \caption{In each subfigure, we plot in the 3D $u$--$v$--$w$ phase space the trajectories for specific locations $x=1,2,...6$ which converge to enclosed orbits.  $\Delta x=0.02$ and $\Delta t=0.05$.}\label{figure3}
\end{figure}

First of all, we take $d_1=5$, $d_2=0.1$, $\mu_1=\mu_2=1$, $\lambda=5$, $\xi=0.1$ and consider (\ref{21}) over interval with length $L=6$, subject to initial condition $(u_0,v_0,w_0)=(\bar u,\bar v,\bar w)+0.001\cos 2 \pi x$.  According to our stability analysis in Proposition \ref{proposition11}, $(\bar u,\bar v,\bar w)$ is unstable if $\chi>\chi_0=\chi^H_{2}\approx 63.2$, and according to our bifurcation analysis and stability results in Theorem \ref{theorem31}, the homogeneous solution loses its stability to time--periodic pattern which has spatial profile $\cos \frac{2\pi x}{L}$; moreover its period is approximately given by $T=\frac{2\pi}{\zeta_0}\approx8$ in (\ref{32}).  In Figure \ref{figure2}, we choose $\chi=80$ and plot $(u,v,w)(x,t)$ for $t\in(0,100)$.  The initial data has a spatial inhomogeneity of the form $\cos 2\pi x$, but the periodic patterns develop according to the spatial profile $\cos \frac{\pi x}{3}$ which is the stable wave mode; moreover the time period of the oscillating patterns matches our theoretical result.  $u$--$v$--$w$ phase spaces are given in Figure \ref{figure3}.

\begin{figure}[h!]
        \centering
\includegraphics[width=\textwidth,height=3in]{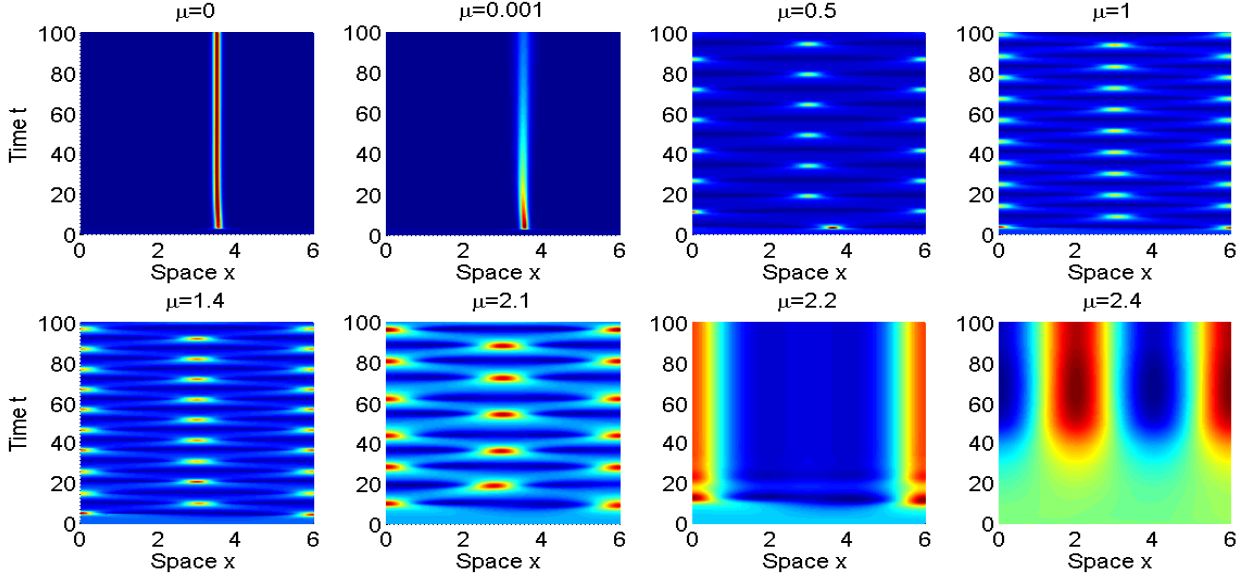}
\caption{Effect of cellular growth on the pattern formation of $u$--species, where we choose $\mu_1=\mu_2$.  System parameters are chosen to be $d_1=8$, $d_2=0.5$, $\chi=130$ and $\xi=0.4$.  Initial data are taken to be small perturbations of $(\bar u,\bar v,\bar w)$.  Space and time grid sizes are $\Delta x=L/500=0.012$ and $\Delta t=0.05$.  We observe that the cellular growth rate $\mu$ supports the formation of periodic patterns.  However, the periodic pattern disappears at $\mu\approx 2.1$, for which we surmise that the oscillating solutions become unstable and develop into a stable stationary pattern.}\label{figure4}
\end{figure}

Figure \ref{figure4} is devoted to illustrating the effects of cellular growth rates $\mu_1$ and $\mu_2$ on the pattern formations in (\ref{21}).  In particular, we choose $\mu_1=\mu_2$ and plot in each subgraph the spatial--temporal behavior of (\ref{21}) as $\mu_i$ increase.

\begin{figure}[h!]
        \centering
\includegraphics[width=\textwidth,height=1.5in]{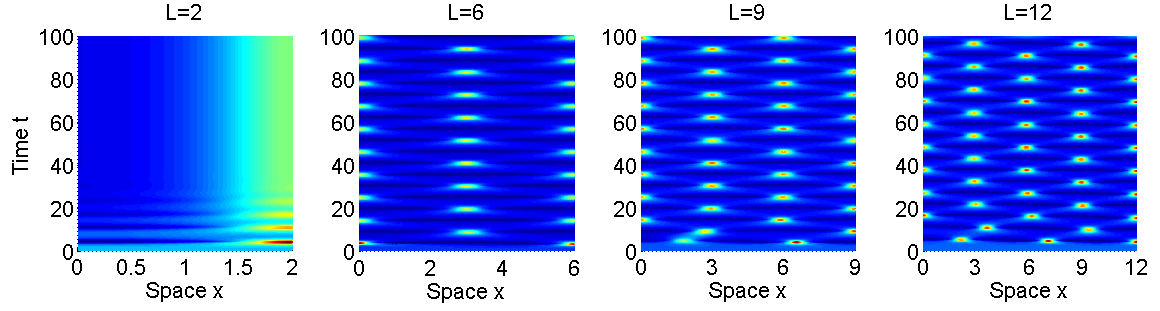}
  \caption{Effect of domain size on the pattern formation of $u$--species.  We choose the system parameters to be the same as those in Figure 3 except that $\chi$ is slightly larger than $\chi_{k_0}$.  $\Delta x=L/500$ and $\Delta t=0.05$ in each graph.  Our simulations support our theoretical findings that large domains support periodic patterns with higher modes, however when the domain size is small, therefore does not exist time--periodic solutions that bifurcate from the homogeneous solution.}\label{figure5}
\end{figure}

If the interval length $L$ is small, we always have that $\chi_k^S<\chi_k^H$, $\forall k\in\mathbb N^+$, then (\ref{21}) does not exist time--periodic solutions through Hopf bifurcation.  To test the effects of $L$ on the pattern formation in (\ref{21}),  Figure \ref{figure5} includes a set of simulations on the spatial--temporal behaviors of solutions to (\ref{21}) over different intervals, subject to the same set of system parameters and initial data as in Figure \ref{figure3}.  We want to point out that when $L$ is small, $(\bar u,\bar v,\bar w)$ is an attractor to (\ref{21}) if $\chi$ is also small; moreover, according to Remark \ref{remark21}, $\chi_0$ approaches infinity as $L$ approaches zero, therefore $\chi$ has to be sufficiently large to support pattern formation when the interval length is small.

\begin{figure}[h!]
        \centering
\includegraphics[width=\textwidth,height=3.5in]{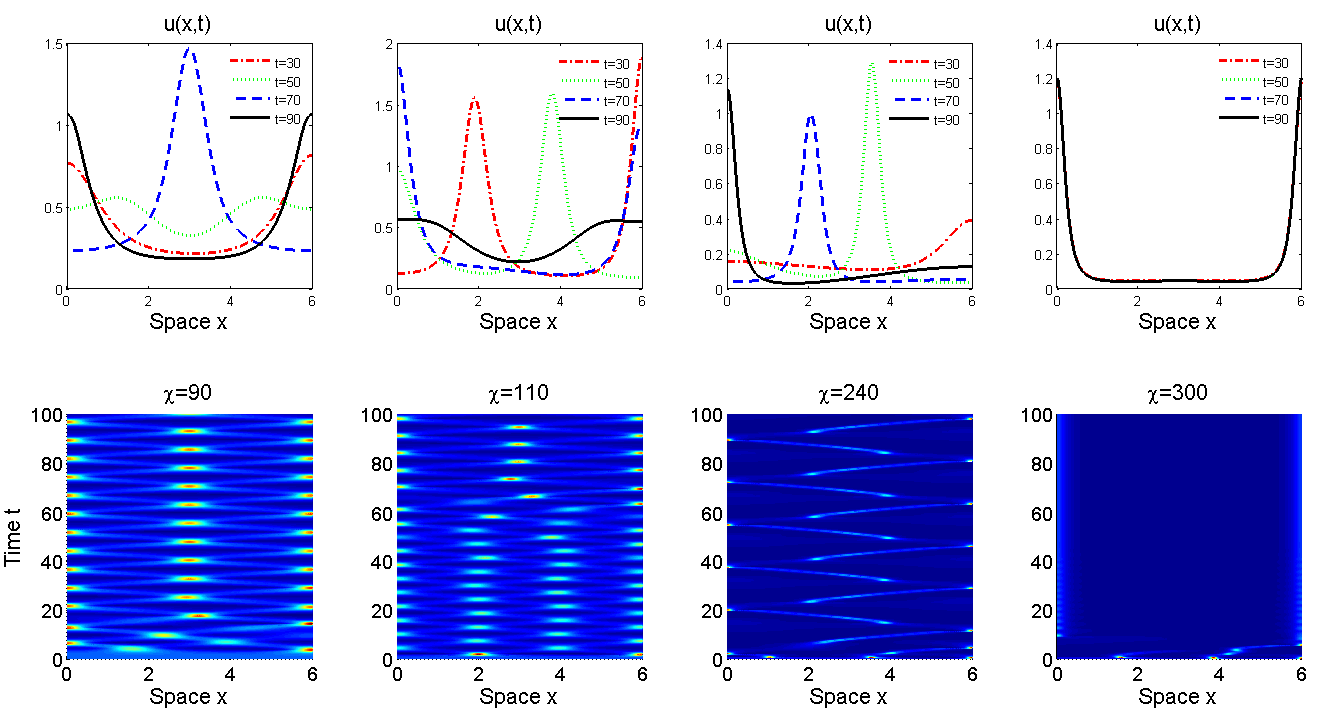}
  \caption{Pattern formation of $u$--species in (\ref{21}) when chemotaxis rate $\chi$ is far away from $\chi_{k_0}$=63.2.}\label{figure6}
\end{figure}

Finally, our numerics in Figure \ref{figure6} are devoted to examining the effects of $\chi$ on the spatial--temporal dynamics of (\ref{21}), when $\chi$ is far away from $\chi_{k_0}=63.2$.  In the plots from left to right, we choose $\chi=90$, 110, 240 and 300 respectively, while the rest parameters and initial data are the same as those in Figure \ref{figure2}.
When $\chi=90$ or $110$, we see that the stable time--periodic solutions have the same profiles $\cos \frac{\pi x}{3}$ as described in Theorem \ref{theorem31}, while a time--periodic spatial pattern with mode $\cos \frac{\pi x}{2}$ is developed for time $t$ up to 60 when $\chi=110$.  We surmise that this oscillating solution is unstable or metastable, and a nonlinear analysis is required to determine its stability.  Moreover, (\ref{21}) has periodic patterns when $\chi=240$ which is far away from $\chi_0\approx 63.2$, and it natural to expect that the existence of this periodic solution is not driven by the linearized instability of the homogeneous solution but the nonlinear cellular growth.

\section{Conclusions and discussions}\label{section6}
In this paper, we study the coupled $3\times3$ system (\ref{11}) which models the spatial--temporal evolution of two competing species and one--attracting chemical.  It has been revealed in \cite{STW,TW} that either the unique positive equilibrium or semi--equilibrium is a global attractor to (\ref{11}) if the chemotaxis coefficients $\chi$ and $\xi$ are small compared to the cellular growth rates $\mu_1$ and $\mu_2$.  Nonconstant positive steady states of (\ref{11}) over $\Omega=(0,L)$ have been studied in \cite{WZYH} by rigorous analysis.  Our results complement the works in \cite{STW,TW,WZYH} by studying its positive time--periodic spatially inhomogeneous solutions observed through the numerical studies in \cite{WZYH}.  Periodic patterns have been experimentally observed in chemotaxis of E. coli or \emph{Dictyostelium discoideum} by many researchers \cite{HP,Ho,Ho2}.  Numerical simulations have been performed to investigate the oscillating patterns by various authors \cite{EIM,PH}, however very few works have been done towards rigorous mathematical analysis of these periodic solutions (the only related work we know is \cite{LSW}).

The starting point of our mathematical analysis of (\ref{11}) is the linearized stability of its positive constant solution $(\bar u,\bar v,\bar w)$, which becomes unstable if $\chi>\chi_0=\min_{k\in\mathbb N^+}\{\chi_k^S,\chi_k^H\}$.  It is proved in \cite{WZYH} that if $\chi_0=\chi^S_{k_1}<\chi_k^H$, $\forall k\in \mathbb N^+$, then the stability of $(\bar u,\bar v,\bar w)$ is lost to nonconstant positive stationary solutions of (\ref{21}) at $\chi^S_{k_1}$ through steady state bifurcation, while all the rest bifurcating solutions around $\chi^S_k$ are unstable if $k\neq k_1$.  The first set of results in this paper state that if $\chi_0=\chi^H_{k_0}<\chi_k^H$, $\chi_k^S$, $\forall k\in \mathbb N^+$, then $(\bar u,\bar v,\bar w)$ loses its stability to time--periodic spatial solutions at $\chi^H_{k_0}$ through Hopf bifurcation, while all the rest Hopf bifurcation branches must be unstable if $k\neq k_0$.  These linearized stability results complete our understanding of the local dynamics of $(\bar u,\bar v,\bar w)$: $(\bar u,\bar v,\bar w)$ is driven unstable by large chemotaxis rate through steady state bifurcation if $\chi_0=\chi^S_{k_1}$ and through Hopf bifurcation if $\chi_0=\chi^H_{k_0}$;  moreover, our stability results provide a complete wave mode selection mechanism for system (\ref{21}) in sense that the only stable bifurcating solution (through Hopf or steady state) must stay on the left--most branch, while all the rest bifurcating branches are academic in that they are all unstable.

Another main achievement of this paper is that we reveal the effect of cellular growth on the dynamics of (\ref{11}) over multi--dimensional bounded domains.  In particular, we showed that the cellular kinetics are responsible for the formation of time--periodic patterns to (\ref{11}), which has no temporal oscillating patterns when $\mu_1=\mu_2=0$.  Our proof is based on the construction of time--monotone Lyapunov functional.  An extra conclusion we have from the Lyapunov functional is that we proved the global existence and boundedness of classical solutions to (\ref{11}) over $\Omega\subset \mathbb R^2$ provided the total cell population is not too large.  Considering system (\ref{11}) over $(0,L)$, we have also studied the effect of the domain size on the pattern formation, which shows that small domain supports stationary patterns while large domain supports time--periodic patterns.  Numerical simulations are implemented to illustrate both stationary and time--periodic patterns to (\ref{11}) which support our theoretical findings.

We know that the critical value $\chi_0=\chi_0(\xi)$ decreases as $\xi$ increases and $\chi_0<0$ if $\xi$ is sufficiently large, therefore only one of $\chi$ and $\xi$ is needed to be large to destabilize $(\bar u,\bar v,\bar w)$.  If $\xi<0$, i.e., species $v$ is repulsive to the chemical gradient, $\chi$ needs to be large to destabilize $(\bar u,\bar v,\bar w)$.  Moreover, the local stability analysis suggests that chemo--attraction destabilizes constant steady states and the chemo--repulsion stabilizes constant steady states.  The constant solution is always stable both when $\chi<0$ and $\xi<0$.  Therefore we surmise that $(\bar u,\bar v,\bar w)$ is also a global attractor of (\ref{11}) if $\chi<0$ and $\xi<0$, while both $\chi$ and $\xi$ are chosen to be positive in \cite{TW}.  This needs an approach totally different from those in this paper.

When $(\bar u,\bar v,\bar w)$ loses its stability at $\chi^H_{k_0}$, a Hopf bifurcation occurs and the stability is lost to a time--periodic solution.  Our analysis of the Hopf bifurcation stability states that it depends on the turning direction of the branch around the equilibrium.  It seems necessary to evaluate $\chi^{H'}_{k_0}(0)$ for this sake which need some very complicated calculations.  The global Hopf bifurcation analysis is also a very interesting problem that worths future exploration.  It is interesting and important to ask what happens when at a later stage the time--periodic solutions lose stability? We refer this to the Poincar\'e map for which \cite{Sat1} is a good reference.

We showed that in the absence of cellular growth, (\ref{11}) does not have any time--periodic solutions since it admits a time--monotone Lyapunov functional.  In light of this, the global existence and boundedness of (\ref{11}) is obtained over two--dimensional domains, where we have assumed the smallness of initial cell populations.  From the viewpoint of mathematical analysis, it is an interesting and important question to study the global existence or blow--ups of (\ref{11}) in higher dimensions when $\mu_1=\mu_2=0$.  In particular, for the global existence on Keller--Segel chemotaxis models, we refer the reader to the very recent survey \cite{BBTW}.  It appears that both $\chi$ and $\xi$ are required positive in the arguments of \cite{TW}, and in light of our stability analysis, we surmise that the solution of (\ref{11}) is always global if $\chi<0$ and $\xi<0$, since chemo--repulsion has smoothing effect like diffusions.  To prove this, one needs an approach totally different from that in \cite{TW}.

\end{document}